%% file: lsmonoids.tex
\documentclass[12pt]{article}
\usepackage{a4, amssymb, latexsym, amsmath, exscale, amscd}
\usepackage[all]{xy}

\newenvironment{evlist}[2]{
\begin{list}{}{
\setlength{\topsep}{0.5ex plus0.2ex minus0.1ex} 
\setlength{\leftmargin}{#1}
\setlength{\itemsep}{#2 plus0.2ex}
\setlength{\parsep}{0ex plus0.2ex} }}
{\end{list}}

\newcommand{\Self}[1]{\mathrm{T}_{#1}}

\newcommand{\Nat}{\mathbb{N}}
\newcommand{\Int}{\mathbb{Z}}
\newcommand{\id}{\mathrm{id}}

\newcommand{\definition}[1]{\textit{#1}}

\newcommand{\proof}{{\textit{Proof}\enspace}}

\newcommand{\eop}{\ \vbox{\hrule
                       \hbox{\vrule
                             \hskip 6pt
                             \vrule height 6pt width 0pt
                             \vrule}%
                       \hrule}%
                     \vspace{\medskipamount}
                }

\newsavebox{\ttt}
\sbox{\ttt}{}
\newcommand{\startsection}[1]
    {\section[#1]{#1}
    \sbox{\ttt}{\thesection\ \ \textsc{#1}}
    \thispagestyle{plain}
}

\newtheorem{lemma}{Lemma}[section]
\newtheorem{proposition}{Proposition}[section]
\newtheorem{theorem}{Theorem}[section]

\newcommand\mythicklines{}
\newcommand\sput[3]{\put(#1,#2){$\scriptstyle{#3}$}}

\newlength{\graphicthick}
\newlength{\graphicmid}
\newlength{\graphicthin}

\begin{document}

\begin{titlepage}

\begin{center}
\phantom{q}
\vspace{60pt}

{\huge Existence of monoids compatible}

\smallskip
{\huge with a family of mappings}

\bigskip\bigskip\bigskip
{\LARGE Chris Preston}

\bigskip
\end{center}

\bigskip
\bigskip
\bigskip
\bigskip
\begin{quote}
These notes, which are a revised version of \cite{preston2}, present an approach to obtaining monoid operations 
which are compatible with a given family of mappings in the sense that the mappings become left translations in 
the monoid. This can be applied to various situations such as the addition on the natural numbers and the integers 
as well as the concatenation of lists. Some of the results can also be found in \cite{preston1}.
\end{quote}

\end{titlepage}

\thispagestyle{empty}

\addtocontents{toc}{\vskip 20pt}
\addtolength{\parskip}{-5pt}
\addtolength{\parskip}{5pt}

\include{lsalg}
\include{tm}
\include{initial}

\end{document}

%% file: lsalg.tex
\startsection{Minimal $\mathcal{L}_S$-algebras}
\label{lsalg}

Let $(X,x_0)$ be a pointed set (i.e., a set $X$ together with a distinguished base-point $x_0 \in X$). 
A binary operation $\bullet$ on $X$ will be called a \definition{monoid operation on $(X,x_0)$}
if $(X,\bullet,x_0)$ is a monoid having $x_0$ as unit element, meaning that $\bullet$ is associative
and $x \bullet x_0 = x_0 \bullet x = x$ for all $x \in X$.

If $\bullet$ is a monoid operation on $(X,x_0)$ then a mapping $u : X \to X$ is said to be a
\definition{(left) translation} in $(X,\bullet,x_0)$ if $u(x) = u(x_0) \bullet x$ for all $x \in X$.
It then it follows from the associativity of $\bullet$ that
\[u(x_1 \bullet x_2) = u(x_0) \bullet (x_1 \bullet x_2) = (u(x_0) \bullet x_1) \bullet x_2
 = u(x_1) \bullet x_2\]
for all $x_1,\,x_2 \in X$, and therefore
\begin{evlist}{28pt}{6pt}
\item[$\mathrm{(\star)}$]
$\ u(x_1) \bullet x_2 = u(x_1 \bullet x_2)$ for all $x_1,\,x_2 \in X$.
\end{evlist}
On the other hand, if $\mathrm{(\star)}$ holds then 
$u(x) = u(x_0 \bullet x) = u(x_0) \bullet x$ for all $x \in X$, and so
$u$ is a translation.

Now let $S$ be a fixed non-empty set. If $X$ is any set then a mapping $f : S \times X \to X$ will also 
be regarded as a family of mappings $\{f_s\}_{s \in S}$, where $f_s : X \to X$ is given by $f_s(x) = f(s,x)$ 
for all $x \in X$. A triple $(X,f,x_0)$ consisting of a non-empty set $X$, a mapping $f : S \times X \to X$  
and an element $x_0 \in X$ will be called an \definition{$\mathcal{L}_S$-algebra}, where the symbol 
$\mathcal{L}$ should be though of as standing for \definition{list}.
The reason for this terminology is that such triples are exactly the algebras -- in the sense of
universal algebra -- associated with the signature for specifying lists of elements from the
set $S$: 
The element $x_0$ corresponds to the empty list and $f_s(x)$ corresponds to the `list' obtained by 
adding the element $s$ to the beginning of the `list' $x$.
(This doesn't mean, though, that the elements in an $\mathcal{L}_S$-algebra have to look anything 
like  real lists.)

The kinds of $\mathcal{L}_S$-algebras we have in mind -- and which are introduced below -- all have the 
additional property of being minimal: An $\mathcal{L}_S$-algebra $\Lambda = (X,f,x_0)$ is said to be 
\definition{minimal} if the only $f$-invariant subset of $X$ containing $x_0$ is $X$ itself, where a subset 
$X'$ of $X$ is \definition{$f$-invariant} if $f_s(X') \subset X'$ for each $s \in S$ (or, what is equivalent, 
if $f(S \times X') \subset X'$). We will thus essentially restrict our attention to such objects.

If $\Lambda = (X,f,x_0)$ is an $\mathcal{L}_S$-algebra then a monoid operation $\bullet$ on $(X,x_0)$ will be 
called \definition{$\Lambda$-compatible} if $f_s$ is a translation in $(X,\bullet,x_0)$ for each $s \in S$.
The aim of these notes is to characterise those minimal $\mathcal{L}_S$-algebras $\Lambda$ for which there 
exists a $\Lambda$-compatible monoid operation. If such an operation exists then it is unique:

\begin{lemma}\label{lemma_lsalg_11}
If $\Lambda$ is minimal then there exists at most one $\Lambda$-compatible monoid operation.
\end{lemma}

\proof 
Let $\Lambda = (X,f,x_0)$ and let $\bullet_1$ and $\bullet_2$ be $\Lambda$-compatible monoid operations.
Then the set $X_0 = \{ x \in X : \mbox{$x \bullet_1 y = x \bullet_2 y$ for all $y \in X$} \}$
contains $x_0$, since $x_0 \bullet_1  y = y = x_0 \bullet_2 y$ for all $y \in X$, and it is $f$-invariant:  
If $x \in X_0$ then for all $y \in X$ it follows from $\mathrm{(\star)}$ that 
$f_s(x) \bullet_1 y = f_s(x \bullet_1 y) = f_s(x \bullet_2 y) = f_s(x) \bullet_2 y$, and thus $f_s(x) \in X_0$ 
for all $s \in S$. Hence $X_0 = X$, since $\Lambda$ is minimal. This shows that ${\bullet_1} = {\bullet_2}$. \eop

Let $\Lambda = (X,f,x_0)$ be a minimal $\mathcal{L}_S$-algebra.
The main result (Theorem~\ref{theorem_lsalg_11}) will present a necessary and sufficient condition for the 
existence of a $\Lambda$-compatible monoid operation in terms of what we call a reflection.
In particular, it follows that if $\Lambda$ is also \definition{commutative},
meaning that $f_s \circ f_t = f_t \circ f_s$ for all $s,\,t \in S$,
then there always exists a $\Lambda$-compatible monoid operation.

If $\Lambda = (X,f,x_0)$ is an $\mathcal{L}_S$-algebra then a mapping $f' : S \times X \to X$ will be called 
a \definition{reflection of $f$ in $\Lambda$} if $f'_s(x_0) = f_s(x_0)$ for all $s \in S$ and 
$f'_s \circ f_t = f_t \circ f'_s$ for all $s,\,t \in S$. If $f'$ is a reflection of $f$ in $\Lambda$ then, 
conversely, $f$ is a reflection of $f'$ in the $\mathcal{L}_S$-algebra $\Lambda' = (X,f',x_0)$.
Of course, $f$ is a reflection of itself in $\Lambda$ if and only if $\Lambda$ is commutative.
There is a related concept for monoid operations: To each monoid operation $\bullet$ on $(X,x_0)$ there
is an associated operation $\bullet'$ given by $x_1 \bullet' x_2 = x_2 \bullet x_1$ for all
$x_1,\,x_2 \in X$, which we also refer to as the \definition{reflection} of $\bullet$.
The relation between $\bullet$ and $\bullet'$ is symmetric in that $\bullet$ is the reflection of
$\bullet'$, and  $\bullet' = \bullet$ if and only if $\bullet$ is commutative.

A reflection is necessary for the existence of a compatible monoid operation, as the next result shows.
Theorem~\ref{theorem_lsalg_11} then states that the converse is also true for minimal $\mathcal{L}_S$-algebras.

\begin{proposition}\label{prop_lsalg_21}
Let $\Lambda = (X,f,x_0)$ be an $\mathcal{L}_S$-algebra and suppose there exists a $\Lambda$-compatible monoid 
operation $\bullet$. Let $f' : S \times X \to X$ be the mapping defined by $f'_s(x) = x \bullet f_s(x_0)$ 
for all $x \in X$, $s \in S$. Then $f'$ is a reflection of $f$ in $\Lambda$, and the reflection 
$\bullet'$ of $\bullet$ is a $\Lambda'$-compatible monoid operation, where $\Lambda' = (X,f',x_0)$.
Moreover, if $\Lambda$ is minimal then so is $\Lambda'$.
\end{proposition}

\proof 
For all $s,\,t \in S$ and all $x \in X$
\begin{eqnarray*}
(f'_s \circ f_t)(x) &=& f'_s(f_t(x)) = f'_s(f_t(x_0) \bullet x) = (f_t(x_0) \bullet x) \bullet f_s(x_0)\\
&=& f_t(x_0) \bullet (x \bullet f_s(x_0)) = f_t(x_0) \bullet f'_s(x) = f_t(f'_s(x)) = (f_t \circ f'_s)(x) \;,
\end{eqnarray*}
i.e., $f'_s \circ f_t = f_t \circ f'_s$. Also $f'_s(x_0) = x_0 \bullet f_s(x_0) = f_s(x_0)$ for all 
$s \in S$, and thus $f'$ is a reflection of $f$ in $\Lambda$.
Moreover, $f'_s(x) = x \bullet f_s(x_0) = f'_s(x_0) \bullet' x$ for all $x \in X$, $s \in S$, and so
$\bullet'$ is a $\Lambda'$-compatible monoid operation.
The proof of the final statement (that if $\Lambda$ is minimal then so is $\Lambda'$), which is not quite
so straightforward, is given later.
\eop

\begin{lemma}\label{lemma_lsalg_21}
If $\Lambda = (X,f,x_0)$ is minimal then there exists at most one reflection of $f$ in $\Lambda$.
\end{lemma}

\proof 
If $f'$ and $f''$ are both reflections of $f$ in $\Lambda$ then the set $X_0$ consisting of those $x \in X$ 
for which $f'_t(x) = f''_t(x)$ for all $t \in S$ contains the element $x_0$, since 
$f'_t(x_0) = f_t(x_0) = f''_t(x_0)$ for all $t \in S$, and it is $f$-invariant: If $x \in X_0$ and $s \in S$ 
then $f'_t(f_s(x)) = f_s(f'_t(x)) =  f_s(f''_t(x)) = f''_t(f_s(x))$ for all $t \in S$ and so $f_s(x) \in X_0$. 
Hence $X_0 = X$, since $\Lambda$ is minimal, which shows that $f'_t = f''_t$ for all $t \in S$, i.e., $f' = f''$.
\eop

\begin{theorem}\label{theorem_lsalg_11}
Let $\Lambda = (X,f,x_0)$ be a minimal $\mathcal{L}_S$-algebra. Then there exists a $\Lambda$-compatible monoid 
operation (which by Lemma~\ref{lemma_lsalg_11} is then unique)
if and only if there is a reflection $f'$ of $f$ in $\Lambda$. 
\end{theorem}

\proof Later.
\eop

Here is the special case of  Theorem~\ref{theorem_lsalg_11} for a commutative $\mathcal{L}_S$-algebra. 

\begin{theorem}\label{theorem_lsalg_21}
Let $\Lambda$ be a minimal commutative $\mathcal{L}_S$-algebra. 
Then there exists a unique $\Lambda$-compatible monoid operation $+$
and this operation is commutative.
\end{theorem}

\proof 
This follows from Theorem~\ref{theorem_lsalg_11} and Proposition~\ref{prop_lsalg_11}~(1).
\eop

Before coming to the examples we give a result which shows the relationship
between properties of a $\Lambda$-compatible monoid operation and properties of the family of mappings 
$\{f_s\}_{s \in S}$.

\begin{proposition}\label{prop_lsalg_11}
Let $\Lambda = (X,f,x_0)$ be a minimal $\mathcal{L}_S$-algebra and suppose there exists a $\Lambda$-compatible 
monoid operation $\bullet$ (which by Lemma~\ref{lemma_lsalg_11} is then unique).
Then:
\begin{evlist}{8pt}{6pt}
\item[(1)]
The monoid $(X,\bullet,x_0)$ is commutative (meaning that $x_1 \bullet x_2 = x_2 \bullet x_1$ for all
\phantom{xxx}$x_1,\,x_2 \in X$) if and only if $\Lambda$ is commutative.

\item[(2)]
The monoid $(X,\bullet,x_0)$ obeys the left cancellation law  (meaning that $x_1 = x_2$ 
\phantom{xxx}whenever $x \bullet x_1 = x \bullet x_2$ 
for some $x \in X$) if and only if $f_s$ is injective for
\phantom{xxx}each $s \in S$.

\item[(3)]
The monoid $(X,\bullet,x_0)$ is a group  if and only if $f_s$ is surjective for each 
\phantom{xxx}$s \in S$, which is the case if and only if $f_s$ is bijective for each $s \in S$.
\end{evlist}
\end{proposition}

\proof Later.
\eop

Note that a monoid $(X,\bullet,x_0)$ obeys the right cancellation law  (meaning that $x_1 = x_2$ whenever 
$x_1 \bullet x = x_2 \bullet x$ for some $x \in X$) if and only if the reflection $\bullet'$ obeys the left
cancellation law. Thus if $\Lambda = (X,f,x_0)$ is a minimal $\mathcal{L}_S$-algebra for which there exists a 
$\Lambda$-compatible monoid operation $\bullet$ then by Theorem~\ref{theorem_lsalg_11} and
Propositions \ref{prop_lsalg_21} and \ref{prop_lsalg_11} 
$(X,\bullet,x_0)$ obeys the right cancellation law if and only if $f'_s$ is injective for each $s \in S$,
where $f'$ is the reflection of $f$ in $\Lambda$.

\bigskip
Let us now look at some typical examples of $\mathcal{L}_S$-algebras.

\bigskip

1.\enskip
Let $\Delta$ be a set consisting of a single element, say $\#$. Consider an $\mathcal{L}_{\Delta}$-algebra 
$\Lambda = (X,f,x_0)$; then $f : \Delta \times X \to X$ can be regarded just as a mapping $f : X \to X$ 
(by identifying $f$ with $f_{\#}$) and $\Lambda$ being minimal means that the only subset $X'$ of $X$ 
containing $x_0$ with $f(X') \subset X'$ is $X$ itself. Here $\Lambda$ is clearly commutative, and so
by Theorem~\ref{theorem_lsalg_21} there exists a unique $\Lambda$-compatible monoid operation $+$ and this
operation is commutative. The most important example here is the $\mathcal{L}_{\Delta}$-algebra 
$(\Nat,\textsf{s},0)$, where $\Nat = \{0,1,\ldots\}$ is the set of natural numbers and 
$\textsf{s} : \Nat \to \Nat$ is the successor operation (with $\textsf{s}(0) = 1$, $\textsf{s}(1) = 2$ and 
so on). The fact that $(\Nat,\textsf{s},0)$ is minimal follows from one of the Peano axioms, namely the axiom 
requiring the \definition{principle of mathematical induction} to hold. The operation $+$ given by 
Theorem~\ref{theorem_lsalg_21} in this case is the addition on $\Nat$: By $\mathrm{(\star)}$ and since $0$ is the 
unit element it follows that
\begin{evlist}{36pt}{6pt}
\item[$\mathrm{(+_0)}$]
$\ 0 + n = n$ for all $n \in \Nat$,
\item[$\mathrm{(+_1)}$]
$\ \textsf{s}(n_1) + n_2 = \textsf{s}(n_1 + n_2)$ for all $n_1,\,n_2 \in \Nat$,
\end{evlist}
and the `equations' $\mathrm{(+_0)}$ and $\mathrm{(+_1)}$ are the usual recursive specification for the addition. 
Proposition~\ref{prop_lsalg_11}~(2) confirms that the cancellation law holds here, since the successor operation 
$\textsf{s} : \Nat \to \Nat$ is injective.

Theorem~\ref{theorem_lsalg_21} shows that the addition on $\Nat$ can be obtained without using the other two 
Peano axioms. These axioms, when stated in terms of an $\mathcal{L}_{\Delta}$-algebra $(X,f,x_0)$ require the 
mapping $f$ to be injective and $f(x) \ne x_0$ to hold for all $x \in X$. What if $(X,f,x_0)$ is minimal but 
one of these axioms does not hold, and so either $x_0 \in f(X)$ or $f$ is not injective? In both cases $X$ is 
finite. If $x_0 \in f(X)$ then $f$ is a bijection and the picture looks like:

\setlength{\graphicthick}{0.1mm}
\setlength{\graphicmid}{0.1mm}
\setlength{\graphicthin}{0.1mm}

\begin{center}
\setlength{\unitlength}{0.8mm}
\begin{picture}(140,50)

\linethickness{\graphicthick}

\linethickness{\graphicthin}

\linethickness{\graphicmid}
\mythicklines

\put(49,24){$\bullet$}
\sput{53}{24}{x_0 = f(x_\ell)}

\sput{58}{0}{x_1 = f(x_0)}
\put(59,4){$\bullet$}

\sput{83}{0}{x_2 = f(x_1)}
\put(84,4){$\bullet$}

\put(59,44){$\bullet$}
\sput{58}{48}{x_\ell}

\put(84,44){$\bullet$}
\put(94,24){$\bullet$}

\put(50,25){\line(1,2){10}}
\put(50,25){\line(1,-2){10}}
\put(60,45){\line(1,0){25}}
\put(60,5){\line(1,0){25}}
\put(85,45){\line(1,-2){10}}
\put(85,5){\line(1,2){10}}

\end{picture}

\end{center}

Here $+$ is really nothing but addition modulo $n$ with $n$ the cardinality of  $X$ and, as confirmed by 
Proposition~\ref{prop_lsalg_11}~(3), in this case the associated monoid is an abelian group. If $f$ is not injective 
then the picture is the following:

\begin{center}
\setlength{\unitlength}{0.8mm}
\begin{picture}(140,55)

\linethickness{\graphicthick}

\linethickness{\graphicthin}

\linethickness{\graphicmid}
\mythicklines

\put(5,25){\line(1,0){90}}

\put(95,25){\line(1,2){10}}
\put(95,25){\line(1,-2){10}}
\put(105,45){\line(1,0){25}}
\put(105,5){\line(1,0){25}}
\put(130,45){\line(1,-2){10}}
\put(130,5){\line(1,2){10}}

\sput{4}{21}{x_0}
\put(4,24){$\bullet$}

\sput{29}{21}{x_1 = f(x_0)}
\put(29,24){$\bullet$}

\sput{73}{21}{x_t}
\put(73,24){$\bullet$}

\sput{98}{24}{\breve{x}_0 = f(\breve{x}_\ell) = f(x_t)}
\put(94,24){$\bullet$}

\sput{103}{0}{\breve{x}_1 = f(\breve{x}_0)}
\put(104,4){$\bullet$}

\sput{103}{48}{\breve{x}_\ell}
\put(104,44){$\bullet$}

\put(129,4){$\bullet$}
\put(129,44){$\bullet$}
\put(139,24){$\bullet$}

\end{picture}

\end{center}

Even if it is surprising that an addition exists in this case it is a simple enough matter to explicitly 
compute what this operation has to be. Since $f$ is not injective Proposition~\ref{prop_lsalg_11}~(2) implies that 
the cancellation law does not hold here.

\bigskip
2.\enskip
Let $(X,x_0)$ be a pointed set, let $f_+ : X \to X$ be a bijection and put $f_- = f_+^{-1}$. Then 
$\Lambda = (X,f,x_0)$ is an $\mathcal{L}_{\pm}$-algebra with ${\pm} = \{+,-\}$, where $f : {\pm} \times X \to X$ 
is given by $f(+,x) = f_+(x)$ and $f(-,x) = f_-(x)$ for all $x \in X$. In particular, $\Lambda$ is commutative.
A very special case of this is the $\mathcal{L}_{\pm}$-algebra $\Lambda_\Int = (\Int,\textsf{s},0)$, where
$\textsf{s}_+(n) = n + 1$ and $\textsf{s}_-(n) = n - 1$ for all $n \in \Int$. It is easily checked that 
$\Lambda_\Int$ is minimal. Thus by Theorem~\ref{theorem_lsalg_21} there is a unique $\Lambda_\Int$-compatible 
monoid operation $+$ which is commutative, and $+$ is uniquely determined by the requirements that 
$m + 0 = 0 = 0 + m$ for all $m \in \Int$ and
\[ \textsf{s}_+(m + n) = \textsf{s}_+(m) + n \mbox{\ \ and\ \ } \textsf{s}_-(m + n) = \textsf{s}_-(m) + n\]
for all $m,\,n \in \Int$. Of course, $+$ is the usual addition on $\Int$. Proposition~\ref{prop_lsalg_11}~(3) 
confirms that $(\Int,+,0)$ is a group, since $\textsf{s}_+$ and $\textsf{s}_-$ are both bijections.

\bigskip
3.\enskip
Let $(M,\bullet,e)$ be a monoid and let $\sigma : M \times M \to M$ be the mapping given by 
$\sigma(a,b) = a \bullet b$. Then $\Lambda = (M,\sigma,e)$ is an $\mathcal{L}_M$-algebra, which is clearly 
minimal since $a = \sigma_a(e)$ for each $a \in M$. Moreover, $\Lambda$ is commutative if and only if the 
monoid $M$ is. This can be generalised somewhat: Again let $(M,\bullet,e)$ be a monoid, let $S$ be a non-empty 
subset of $M$ and let $\sigma : S \times M \to M$ be the mapping given by $\sigma(a,b) = a \bullet b$.
Then $\Lambda = (M,\sigma,e)$ is an $\mathcal{L}_S$-algebra, and it is easy to see that $\Lambda$ is minimal 
if and only if the only submonoid of $M$ containing $S$ is $M$ itself. Moreover, $\Lambda$ is commutative if 
and only if $a \bullet b = b \bullet a$ for all $a,\,b \in S$ (and if $\Lambda$ is minimal then this is the 
case if and only if $M$ is commutative). Here the operation $\bullet$ itself is clearly a $\Lambda$-compatible 
monoid operation, since $\sigma_a(b) = a \bullet b = (a \bullet e) \bullet b = \sigma_a(e) \bullet b$
for all $a \in S$, $b \in M$. Moreover, the reflection $\sigma'$ of $\sigma$ in $\Lambda$ is given by 
$\sigma'(a,b) = b \bullet a$.

\bigskip
4.\enskip
Denote by $S^*$ the set of all finite lists of elements from $S$, let $\varepsilon$ be the empty list and 
$\triangleleft : S \times S^* \to S^*$ be the mapping such that $\triangleleft_s$ is the operation of adding 
the element $s$ to the beginning of a list. More precisely, $S^* = \bigcup_{n \ge 0} S^n$, with the element
$(s_1,\ldots,s_n)$ of $S^n$ usually written as $s_1\ \cdots\ s_n$, $\varepsilon$ is the single element in $S^0$ 
and $\triangleleft_s(s_1\ \cdots\ s_n) = s\ s_1\ \cdots\ s_n$, with 
$\triangleleft_s(\varepsilon) = s \in S^1 \subset S^*$. Then $\Lambda = (S^*,\triangleleft,\varepsilon)$ is 
the eponymous $\mathcal{L}_S$-algebra and it is easy to see that $\Lambda$ is minimal. Note that if $S$ 
consists of more than one element then $\Lambda$ is not commutative. Here the concatenation operation 
$\bullet$ given by
\[s_1\ \cdots\ s_m \,\bullet\, s'_1\ \cdots\ s'_n = s_1\ \cdots\ s_m \ s'_1\ \cdots\ s'_n\]
is a $\Lambda$-compatible monoid operation. Moreover, the reflection $\triangleleft'$ of $\triangleleft$ in 
$\Lambda$ is the mapping such that $\triangleleft'_s$ is the operation of adding the element $s$ to the end 
of a list, i.e., $\triangleleft'_s(s_1\ \cdots\ s_n) = s_1\ \cdots\ s_n\ s$. We will deal with this example in 
more detail in Section~\ref{initial}.

\bigskip
5.\enskip
Let $n \ge 1$ and define a mapping $f : S \times S^n \to S^n$ by
\[   f(s,(s_1,\dots,s_n)) = (s,s_1,\ldots,s_{n-1})\;; \]
let $x_0$ be any element of $S^n$. Then $\Lambda = (S^n,f,x_0)$ is an $\mathcal{L}_S$-algebra, which is minimal 
since $(s_1,\ldots,s_n) = f_{s_1}(f_{s_2}( \cdots f_{s_n}(x)\cdots))$ for all $(s_1,\ldots,s_n) \in S^n$, 
$x \in S^n$, and in particular with $x = x_0$. If $S$ consists of more than one element then $\Lambda$ is not 
commutative, and in this case there is no $\Lambda$-compatible monoid operation: Suppose $\bullet$ were 
$\Lambda$-compatible; then for all $x_1,\,x_2 \in S^n$ with $x_1 = (s_1,\ldots,s_n)$ it would follow from 
$\mathrm{(\star)}$ that
\begin{eqnarray*} x_1 \bullet x_2 &=& (s_1,\ldots,s_n) \bullet x_2 
= f_{s_1}(f_{s_2}( \cdots f_{s_n}(x_0)\cdots)) \bullet x_2\\
&=& f_{s_1}(f_{s_2}( \cdots f_{s_n}(x_0)\cdots)\bullet x_2)
= \cdots = f_{s_1}(f_{s_2}( \cdots f_{s_n}(x_2)\cdots)) = x_1
\end{eqnarray*}
and in particular that $x = x_0 \bullet x = x_0$ for all $x \in S^n$.

\bigskip
6.\enskip
An alternative description of $\mathcal{L}_S$-algebras is that they are semiautomata with input alphabet $S$ 
and a specified initial state: By definition a \definition{semiautomaton} is a triple $(Q,\Sigma,\delta)$ 
consisting of a set $Q$ (the set of states), a set $\Sigma$ (the input alphabet, which is usually finite) and 
a mapping $\delta : Q \times \Sigma \to Q$ (the transition function). If $(Q,\Sigma,\delta)$ is a 
semiautomaton and $q_0 \in Q$, which can be considered as an initial state, then $(Q,\delta',q_0)$ is an 
$\mathcal{L}_\Sigma$-algebra, where $\delta' : \Sigma \times Q \to Q$ is obtained by transposing the arguments 
of $\delta$, i.e., $\delta'(\sigma,q) = \delta(q,\sigma)$ for all $\sigma \in \Sigma$, $q \in Q$.
Conversely, if $(X,f,x_0)$ is an $\mathcal{L}_S$-algebra then $(X,S,f')$ is a semiautomaton with input alphabet 
$S$, where again $f'$ is obtained by transposing the arguments of $f$, and $x_0$ is an initial state.

\bigskip

We now give the proofs which were omitted above. In Section~\ref{tm} we give an alternative approach to proving 
these results. This is based on Cayley's theorem (in its version for monoids). 

\medskip

\textit{Proof of the final statement in Proposition~\ref{prop_lsalg_21}:}\enskip
We are assuming $\Lambda$ is minimal and  must show that $\Lambda'$ is also minimal.
Let $X'$ be the least $f'$-invariant subset of $X$ containing $x_0$.
We first show that $X'$ is a submonoid of $(X,\bullet,x_0)$, i.e., $y \bullet x \in X'$ for all $y,\,x \in X'$.
Let $X_0 = \{ x \in X' : \mbox{$x \bullet' y \in X'$ for all $y \in X'$} \}$, and so in particular
$x_0 \in X_0$, since $x_0 \bullet' y = y \in X'$ for all $y \in X'$. Let $x \in X_0$ and $s \in S$; 
if $y \in X'$ then $x \bullet' y \in X'$ and therefore $f'_s(x) \bullet' y = f'_s(x \bullet' y) \in X'$, since 
$X'$ is $f'$-invariant, i.e., $f'_s(x) \in X'$. 
Hence $X_0$ is an $f'$-invariant subset of $X'$ 
containing
$x_0$, which implies $X_0 = X'$, since $X'$ is the least $f'$-invariant subset of $X$ containing $x_0$. 
This establishes that
$y \bullet x = x \bullet' y \in X'$ for all $x,\,y \in X'$. 
Next consider $x \in X'$ and $s \in S$; then $f'_s(x_0) \in X'$, since $X'$ is $f'$-invariant and contains $x_0$, 
and $f_s(x) = f_s(x_0) \bullet x = f'_s(x_0) \bullet x \in X'$. Hence $X'$ is $f$-invariant and contains
$x_0$, which implies that $X' = X$, since $\Lambda$ is minimal. Therefore $\Lambda'$ is also minimal. 
\eop

We now prepare for the proof of Theorem~\ref{theorem_lsalg_11}.
In what follows let $\Lambda = (X,f,x_0)$ be a minimal $\mathcal{L}_S$-algebra.

\begin{lemma}\label{lemma_lsalg_31}
A binary operation  $\,\bullet$ on $X$ satisfying 
\begin{evlist}{12pt}{6pt}
\item[$\mathrm{(\bullet_0)}$]
$\ x_0 \bullet x = x$ for all $x \in X$,
\item[$\mathrm{(\bullet_1)}$]
$\ f_s(x_1) \bullet x_2 = f_s(x_1 \bullet x_2)$ for all $x_1,\,x_2 \in X$ and all $s \in S$
\end{evlist}
is a $\Lambda$-compatible monoid operation on $(X,x_0)$.
\end{lemma}

\proof 
Let $\bullet$ be a binary operation satisfying $\mathrm{(\bullet_0)}$ and $\mathrm{(\bullet_1)}$. Then 
$\bullet$ is associative: The set
$X_0 = \{ x \in X : \mbox{$x \bullet (x_1 \bullet x_2) = (x \bullet x_1) \bullet x_2$
                    for all $x_1,\,x_2 \in X$} \}$
contains $x_0$, since $\mathrm{(\bullet_0)}$ implies 
$x_0 \bullet (x_1 \bullet x_2) = x_1 \bullet x_2 = (x_0 \bullet x_1) \bullet x_2$ for all $x_1,\,x_2 \in X$, 
and it is $f$-invariant: If $x \in X_0$ then $x \bullet (x_1 \bullet x_2) = (x \bullet x_1) \bullet x_2$ for 
all $x_1,\,x_2 \in X$ and therefore by $\mathrm{(\bullet_1)}$ 
\begin{eqnarray*}
f_s(x) \bullet (x_1 \bullet x_2) &=& f_s(x \bullet (x_1 \bullet x_2))\\ 
&=& f_s((x \bullet x_1) \bullet x_2) 
= f_s(x \bullet x_1) \bullet x_2 = (f_s(x) \bullet x_1) \bullet x_2 
\end{eqnarray*}
and thus $f_s(x) \in X_0$ for all $s \in S$. Hence $X_0 = X$, since $\Lambda$ is minimal. This shows 
$\bullet$ is associative.

Similarly, the set $X_0 = \{ x \in X : x \bullet x_0 = x \}$ contains $x_0$, since
$x_0 \bullet  x_0 = x_0$ by $\mathrm{(\bullet_0)}$, and it is $f$-invariant: If $x \in X_0$ then 
$f_s(x) \bullet x_0 = f_s(x \bullet x_0) = f_s(x)$ by $\mathrm{(\bullet_1)}$, and thus $f_s(x) \in X_0$ for all 
$s \in S$. Hence $X_0 = X$, since $\Lambda$ is minimal. This shows that $x \bullet x_0 = x$ for all $x \in X$,
which in turn implies that $x \bullet x_0 = x_0 \bullet x = x$ for all $x \in X$, since $\mathrm{(\bullet_0)}$ 
holds, and so we have established that $\bullet$ is a monoid operation on $(X,x_0)$.

Finally, $\bullet$ is $\Lambda$-compatible because $\mathrm{(\bullet_1)}$ is the same as $\mathrm{(\star)}$ 
and $\mathrm{(\star)}$ holding for a mapping is equivalent to it being a translation. \eop

If $x \in X$ then a mapping  $\varrho : X \to X$ will be called \definition{$x$-allowable} if 
$\varrho(x_0) = x$ and $f_s \circ \varrho = \varrho \circ f_s$ for all $s \in S$.

\begin{lemma}\label{lemma_lsalg_41}
If there exists a reflection $f'$ of $f$ in $\Lambda$ then for each $x \in X$ there exists a unique 
$x$-allowable mapping $\varrho_x : X \to X$.
\end{lemma}

\proof 
Consider the set $X_0$ consisting of those elements $x \in X$ for which there exists an $x$-allowable mapping. 
Then $x_0 \in X_0$, since $\id_X$ is $x_0$-allowable, and $X_0$ is $f$-invariant: Let $x \in X_0$ with 
$x$-allowable mapping $\varrho : X \to X$, let $s \in S$ and put $\varrho' = \varrho \circ f'_s$. Then
$\varrho'(x_0) = \varrho(f'_s(x_0)) = \varrho(f_s(x_0)) = f_s(\varrho(x_0)) = f_s(x)$ and for all $t \in S$
\begin{eqnarray*}
f_t \circ \varrho' &=& f_t \circ (\varrho \circ f'_s) = (f_t \circ \varrho) \circ f'_s \\
&=& (\varrho \circ f_t) \circ f'_s = \varrho \circ (f_t \circ  f'_s) = \varrho \circ (f'_s \circ f_t) 
= (\varrho \circ f'_s) \circ f_t = \varrho' \circ f_t
\end{eqnarray*}
and so $\varrho'$ is $f_s(x)$-allowable, i.e., $f_s(x) \in X_0$. Thus  $X_0 = X$, since $\Lambda$ is minimal. 
This shows that for each $x \in X$ there exists a mapping $\varrho : X \to X$ with $\varrho(x_0) = x$ and 
$f_t \circ \varrho = \varrho \circ f_t$ for all $t \in S$.

Now for the uniqueness: Let $x' \in X$ and $\varrho_1,\,\varrho_2$ be $x'$-allowable mappings. Then 
$X_0 = \{ x \in X : \varrho_1(x) = \varrho_2(x) \}$ contains $x_0$, since $\varrho_1(x_0) = x' = \varrho_2(x_0)$,
and it is $f$-invariant, since if $x \in X_0$ then 
$\varrho_1(f_s(x)) = f'_s(\varrho_1(x)) = f'_s(\varrho_2(x)) = \varrho_2(f_s(x))$ and so $f_s(x) \in X_0$ for 
all $s \in S$. Again this implies that $X_0 = X$, which shows that $\varrho_1 = \varrho_2$. \eop

\textit{Proof of Theorem~\ref{theorem_lsalg_11}:}\enskip
We are assuming there exists a reflection $f'$ of $f$ in $\Lambda$, and so by Lemma~\ref{lemma_lsalg_41} there exists 
for each $x \in X$ a unique $x$-allowable mapping $\varrho_x : X \to X$. Define a binary operation $\bullet$ 
on $X$ by letting $x' \bullet x = \varrho_x(x')$ for all $x,\,x' \in X$. Then 
$x_0 \bullet x = \varrho_x(x_0) = x$ for all $x \in X$ and 
\[f_s(x_1) \bullet x_2 = \varrho_{x_2}(f_s(x_1)) = f_s(\varrho_{x_2}(x_1)) = f_s(x_1 \bullet x_2)\] 
for all $x_1,\,x_2 \in X$ and all $s \in S$, and hence $\mathrm{(\bullet_0)}$ and $\mathrm{(\bullet_1)}$ hold. Thus 
by Lemma~\ref{lemma_lsalg_31} $\bullet$ is a $\Lambda$-compatible monoid operation on $(X,x_0)$. 
Conversely, if there exists a $\Lambda$-compatible monoid operation on $(X,x_0)$ then 
Proposition~\ref{prop_lsalg_21} shows that there exists a reflection of $f$ in $\Lambda$.
\eop

\textit{Proof of  Proposition~\ref{prop_lsalg_11}:}\enskip
(1)\enskip
Suppose  first that $(X,\bullet,x_0)$ is commutative and let $s,\,t \in S$. Then
\begin{eqnarray*}
(f_s \circ f_t)(x) &=& f_s(f_t(x)) = f_s(x_0) \bullet f_t(x)  = f_s(x_0) \bullet f_t(x_0) \bullet x\\
&=& f_t(x_0) \bullet f_s(x_0) \bullet x = f_t(x_0) \bullet f_s(x) = f_t(f_s(x)) = (f_t \circ f_s)(x) 
\end{eqnarray*}
for all $x \in X$ and hence $f_s \circ f_t = f_t \circ f_s$, which shows that $\Lambda$ is commutative.
Suppose conversely $\Lambda$ is commutative.
We first show that $f_s(y) \bullet x = y \bullet f_s(x)$ for all $x,\,y \in X$, $s \in S$, and for this
fix $x \in X$ and $s \in S$ and consider the set $X' = \{ y \in X : f_s(y) \bullet x = y \bullet f_s(x) \}$.
Then $f_s(x_0) \bullet x = f_s(x) = x_0 \bullet f_s(x)$, and so $x_0 \in X'$, and if $y \in X'$ and $t \in S$ then
\[ f_s(f_t(y)) \bullet x = f_t(f_s(y)) \bullet x = f_t(f_s(y) \bullet x) = f_t(y \bullet f_s(x))
= f_t(y) \bullet f_s(x) \]
and so $f_t(y) \in X'$. Thus $X'$ is $f$-invariant and contains $x_0$ and therefore $X' = X$, since $\Lambda$ is 
minimal, which means that $f_s(y) \bullet x = y \bullet f_s(x)$ for all $x,\,y \in X$, $s \in S$.
Now consider the set $X_0 = \{ x \in X : \mbox{$x \bullet y  = y \bullet x$ for all $y \in X$} \}$ and so in 
particular $x_0 \in X_0$. If $x \in X_0$ and $s \in S$ then
\[ f_s(x) \bullet y = f_s(x \bullet y) = f_s(y \bullet x) = f_s(y) \bullet x = y \bullet f_s(x) \]
for all $y \in X$, i.e., $f_s(x) \in X_0$. Hence
$X_0$ is $f$-invariant and contains $x_0$ and therefore $X_0 = X$, since $\Lambda$ is 
minimal. This shows $(X,\bullet,x_0)$ is commutative.

(2)\enskip
Suppose first that $f_s$ is injective for each $s \in S$ and let
\[ X_0 = \{ x \in X : \mbox{$x_1 = x_2$ whenever $x \bullet x_1 = x \bullet x_2$} \}\;;\]
in particular $x_0 \in X_0$. Consider $x \in X_0$ and $s \in S$; if $f_s(x) \bullet x_1 = f_s(x) \bullet x_2$ 
then $f_s(x \bullet x_1) = f_s(x) \bullet x_1 = f_s(x) \bullet x_2 = f_s(x \bullet x_2)$, hence
$x \bullet x_1 = x \bullet x_2$, since $f_s$ is injective, and so $x_1 = x_2$, since $x \in X_0$, i.e.,
$f_s(x) \in X_0$. Thus $X_0$ is $f$-invariant and contains $x_0$ and therefore $X_0 = X$, since $\Lambda$ is 
minimal. This shows that $(X,\bullet,x_0)$ obeys the left cancellation law. Suppose conversely that 
$(X,\bullet,x_0)$ does obey the left cancellation law, let $s \in S$ and let $x_1,\,x_2 \in X$ with 
$f_s(x_1) = f_s(x_2)$. Then $f_s(x_0) \bullet x_1 = f_s(x_1) = f_s(x_2) = f_s(x_0) \bullet x_2$ and therefore 
$x_1 = x_2$, which implies that $f_s$ is injective.

(3)\enskip
Suppose first that $f_s$ is surjective for each $s \in S$. Let
\[ X_0 = \{ x \in X : \mbox{for each $y \in X$ there exists $z \in X$ such that $x \bullet z = y$} \}\;,\]
and so in particular $x_0 \in X_0$, since $x_0 \bullet y = y$ for all $y \in X$.
Consider $x \in X_0$ and $s \in S$, and let $y \in X$; since $f_s$ is surjective there exists $y' \in X$
such that $f_s(y') = y$ and since $x \in X_0$ there then exists $z \in X$ with $x \bullet z = y'$.
Hence $f_s(x) \bullet z = f_s(x \bullet z) = f_s(y') = y$, which shows that $f_s(x) \in X_0$.
Therefore $X_0$ is $f$-invariant and contains $x_0$ and thus $X_0 = X$, since $\Lambda$ is 
minimal. In particular, for each $x \in X$ there exists $y \in X$ such that $x \bullet y = x_0$, and 
this implies that $(X,\bullet,x_0)$ is a group. (Let $x \in X$; then there exists $x' \in X$
with $x \bullet x' = x_0$ and $x'' \in X$ with $x' \bullet x'' = x_0$ and so
$x = x \bullet x_0 = x \bullet x' \bullet x'' = x_0 \bullet x'' = x''$, i.e.,
$x' \bullet x = x_0$.) Suppose conversely that $(X,\bullet,x_0)$ is a group, let $s \in S$ and let $x \in X$;
then there exists $y \in X$ such that $f_s(x_0) \bullet y = x$, i.e., such that
$f_s(y) = f_s(x_0) \bullet y = x$. Hence $f_s$ is surjective. 

Finally, if $f_s$ is surjective for each $s \in S$ then $(X,\bullet,x_0)$ is a group, and 
the left cancellation law holds in any group. Thus by (2) $f_s$ is also injective for each $s \in S$.
\eop


%% file: tm.tex
\startsection{Transformation monoids}
\label{tm}

In what follows let $(X,x_0)$ be a fixed pointed set. Denote the set of all mappings of $X$ into itself by 
$\Self{X}$; we thus have the monoid $(\Self{X},\circ,\id_X)$, where $\circ$ is functional composition and 
$\id_X$ is the identity mapping. The submonoids of $\,\Self{X}$ are often referred to as
\definition{transformation monoids}.

The results of Section~\ref{lsalg} will be established using properties of certain of these submonoids.
To be a bit more definite: For each mapping $f : S \times X \to X$ let $M_f$ be the least submonoid of 
$\Self{X}$ containing $f_s$ for each $s \in S$ (i.e., $M_f$ is the intersection of all such submonoids).
Then the statements about the existence of a $\Lambda$-compatible monoid operation for a minimal 
$\mathcal{L}_S$-algebra $\Lambda = (X,f,x_0)$ can all be deduced from properties of $M_f$ and related submonoids.

Recall that a binary operation $\bullet$ on $X$ is a monoid operation on $(X,x_0)$ if $(X,\bullet,x_0)$ is a 
monoid having $x_0$ as unit element.
Moreover, if $\bullet$ is a monoid operation on $(X,x_0)$ then a mapping $u \in \Self{X}$ is 
a translation in $(X,\bullet,x_0)$ if $u(x) = u(x_0) \bullet x$ for all $x \in X$, and in this case it follows 
from the associativity of $\bullet$ that
\begin{evlist}{28pt}{6pt}
\item[$\mathrm{(\star)}$]
$\ u(x_1) \bullet x_2 = u(x_1 \bullet x_2)$ for all $x_1,\,x_2 \in X$.
\end{evlist}

Let $\Phi : \Self{X} \to X$ be the evaluation mapping at $x_0$ given by $\Phi(u) = u(x_0)$ for each 
$u \in \Self{X}$. The restriction of this mapping to a subset $A$ of $\Self{X}$ will be denoted by 
$\Phi_A$; in particular there is then the mapping $\Phi_M : M \to X$ for each submonoid $M$ of $\Self{X}$.

We start with Cayley's theorem (in its version for monoids); this provides the key to the approach we are going to 
take here. For each monoid operation $\bullet$ on $(X,x_0)$ define a mapping $\Psi_\bullet : X \to \Self{X}$ by 
letting
\[ \Psi_\bullet(x)(x') = x \bullet x'\] 
for all $x,\,x' \in X$, and put $M_\bullet = \Psi_\bullet(X)$.

\begin{theorem}[Cayley's theorem]\label{theorem_tm_cayley}
Let $\bullet$ be a monoid operation on $(X,x_0)$; then the following hold:
\begin{evlist}{8pt}{6pt}
\item[(1)]
$\Psi_\bullet$ is an injective homomorphism from 
$(X,\bullet,x_0)$ to $(\Self{X},\circ,\id_X)$; thus $M_\bullet$ is 
\phantom{xxx}a submonoid of $\,\Self{X}$ and  $\Psi_\bullet : (X,\bullet,x_0) \to (M_\bullet,\circ,\id_X)$
is an isomorphism. 

\item[(2)]
$M_\bullet$ consists exactly of the translations in $(X,\bullet,x_0)$, and so in particular
\phantom{xxx}the set of translations is a submonoid of $\,\Self{X}$.

\item[(3)]
The inverse of the isomorphism $\Psi_\bullet$ is the mapping $\Phi_{M_\bullet} : M_\bullet \to X$, and
\phantom{xxx}hence $\Phi_{M_\bullet} : (M_\bullet,\circ,\id_X) \to (X,\bullet,x_0)$ is an isomorphism.
\end{evlist}
\end{theorem}

\proof
(1)\enskip The mapping $\Psi_\bullet$ is a homomorphism since if $x_1,\,x_2 \in X$ then 
\begin{eqnarray*}
\Psi_\bullet(x_1 \bullet x_2)(x) = (x_1 \bullet x_2) \bullet x 
&=& x_1 \bullet (x_2 \bullet x) = \Psi_\bullet(x_1)(x_2 \bullet x) \\
&=& \Psi_\bullet(x_1)(\Psi_\bullet(x_2)(x)) = (\Psi_\bullet(x_1) \circ \Psi_\bullet(x_2))(x) 
\end{eqnarray*}
for all $x \in X$, i.e., $\Psi_\bullet(x_1 \bullet x_2) = \Psi_\bullet(x_1) \circ \Psi_\bullet(x_2)$, and 
$\Psi_\bullet(x_0)(x) = x_0 \bullet x = x$ for all $x \in X$, i.e., $\Psi_\bullet(x_0) = \id_X$. It is injective, 
since if $\Psi_\bullet(x_1) = \Psi_\bullet(x_2)$ then
\[ x_1 = x_1\bullet x_0 = \Psi_\bullet(x_1)(x_0) = \Psi_\bullet(x_2)(x_0) = x_2 \bullet x_0 = x_2\;. \]

(2)\enskip If $u \in M_\bullet$ then $u = \Psi_\bullet(y)$ for some $y \in X$ and thus
\[  u(x) = \Psi_\bullet(y)(x) = y \bullet x = (y \bullet x_0) \bullet x
     = \Psi_\bullet(y)(x_0)\, \bullet x = u(x_0) \bullet x \]
for all $x \in X$, i.e., $u$ is a translation in $(X,\bullet,x_0)$. Conversely, suppose $u \in \Self{X}$ is a 
translation in $(X,\bullet,x_0)$ and put $v = \Psi_\bullet(u(x_0))$; then $v \in M_\bullet$ and
\[  v(x) = \Psi_\bullet(u(x_0))(x) = u(x_0) \bullet x = u(x) \]
for all $x \in X$, i.e., $v = u$ and so $u \in M_\bullet$. This shows $M_\bullet$ is the set of translations 
in $(X,\bullet,x_0)$.

(3)\enskip For each $x \in X$ we have
$(\Phi_{M_\bullet} \circ \Psi_\bullet)(x) = \Phi_{M_\bullet}(\Psi_\bullet(x)) = x \bullet x_0 = x = \id_X(x)$,
therefore $\Phi_{M_\bullet} : M_\bullet \to X$ is the set-theoretic inverse of $\Psi_\bullet$ and hence also the 
inverse of the monoid isomorphism. \eop

\begin{lemma}\label{lemma_tm_11}
If $\,\bullet_1$ and $\,\bullet_2$ are monoid operations on $(X,x_0)$ with $M_{\bullet_1} = M_{\bullet_2}$ then 
$\bullet_1 = \bullet_2$.
\end{lemma}

\proof 
Put $M = M_{\bullet_1} = M_{\bullet_2}$.
By Theorem~\ref{theorem_tm_cayley}~(3) $\Phi_M$ is the inverse of both $\Psi_{\bullet_1}$ and $\Psi_{\bullet_2}$
and therefore $\Psi_{\bullet_1} = \Psi_{\bullet_2}$ (considered as mappings from $X$ to $M$). It thus follows that
$x \bullet_1 x' = \Psi_{\bullet_1}(x)(x') = \Psi_{\bullet_2}(x)(x') = x \bullet_2 x'$
for all $x,\,x' \in X$, i.e., $\bullet_1 = \bullet_2$.
\eop

Theorem~\ref{theorem_tm_cayley} 
implies that if $\Lambda = (X,f,x_0)$ is an $\mathcal{L}_S$-algebra for which there exists a 
$\Lambda$-compatible monoid operation $\bullet$ then $f_s \in M_\bullet$ for each $s \in S$ and therefore
$M_f \subset M_\bullet$.

Let us say that a subset $A$ of $\Self{X}$ is \definition{$x_0$-minimal} if the only $A$-invariant subset of 
$X$ containing $x_0$ is $X$ itself, where $X'$ is \definition{$A$-invariant} if it is $u$-invariant (i.e., 
$u(X') \subset X'$) for each $u \in A$. Now the set $\{ u \in \Self{X} : \mbox{$X'$ is $u$-invariant}\, \}$ 
is clearly a submonoid of $\,\Self{X}$ for each $X' \subset X$; thus if $\Lambda = (X,f,x_0)$ is an 
$\mathcal{L}_S$-algebra then a 
subset $X'$ of $X$ is $f$-invariant if and only if it is $M_f$-invariant. Hence $\Lambda$ is minimal if and 
only if the submonoid $M_f$ is $x_0$-minimal.

\begin{lemma}\label{lemma_tm_21}
A submonoid $M$ of $\,\Self{X}$ is $x_0$-minimal if and only if the mapping $\Phi_M$ is surjective.
\end{lemma}

\proof 
Put $X_0 = \Phi_M(M)$; then $x_0 = \Phi_M(\id_X) \in X_0$, and if $x = \Phi_M(v) \in X_0$ then 
$u(x) = u(v(x_0)) = \Phi_M(u \circ v) \in X_0$ for all $u \in M$; hence $X_0$ is an $M$-invariant subset 
of $X$ containing $x_0$. But each element of $X_0$ has the form $v(x_0)$ for some $v \in M$ and so lies in any 
any $M$-invariant subset of $X$ containing $x_0$, and this implies $X_0$ is the least $M$-invariant subset of 
$X$ containing $x_0$. Hence $X_0 = X$ (i.e., $\Phi_M(M) = X$) if and only if $M$ is $x_0$-minimal, and 
$\Phi_M(M) = X$ is the same as $\Phi_M$ being surjective.
\eop

For each subset $A$ of $\Self{X}$ denote by $Z_A$ the \definition{centraliser of $A$ in $\Self{X}$}, i.e.,
\[Z_A = \{ u \in \Self{X} : \mbox{$u \circ v = v \circ u$ for all $v \in A$} \}\;.\]
The centraliser $Z_A$ is a submonoid of $\Self{X}$, since $\id_X \circ u = u = u \circ \id_X$ for all $u \in A$
and if $v_1,\,v_2 \in Z_A$ then $v_1 \circ v_2 \circ u = v_1 \circ u \circ v_2 = u \circ v_1 \circ v_2$ for all 
$u \in A$.

\begin{lemma}\label{lemma_tm_31}
If $M$ is an $x_0$-minimal submonoid of $\,\Self{X}$ then the mapping $\Phi_{Z_M}$ is injective.
\end{lemma}

\proof 
Let $u_1,\,u_2 \in Z_M$ with $\Phi_{Z_M}(u_1) = \Phi_{Z_M}(u_2)$, i.e., with $u_1(x_0) = u_2(x_0)$. 
Then the set $X_0 = \{ x \in X : u_1(x) = u_2(x) \}$ contains $x_0$, and it is $M$-invariant, since if 
$u_1(x) = u_2(x)$ then $u_1(v(x)) = v(u_1(x)) = v(u_2(x)) = u_2(v(x))$ for all $v \in M$. Thus $X_0 = X$, 
since $M$ is $x_0$-minimal, i.e., $u_1 = u_2$, which implies that $\Phi_{Z_M}$ is injective. 
\eop

A subset $A$ of a submonoid $M$ of $\Self{X}$ is called a \definition{generator of $M$} if $M$ 
is the least submonoid of $\Self{X}$ containing $A$. In particular, if $\Lambda = (X,f,x_0)$ is an 
$\mathcal{L}_S$-algebra and $f_S = \{ u \in \Self{X} : \mbox{$u = f_s$ for some $s \in S$} \}$, then $f_S$ is 
a generator of $M_f$.
Recall that to each monoid operation $\bullet$ on $(X,x_0)$ the reflection 
$\bullet'$ of $\bullet$ is the monoid operation given by $x_1 \bullet' x_2 = x_2 \bullet x_1$ for all 
$x_1,\,x_2 \in X$.

\begin{theorem}\label{theorem_tm_21}
Let $M$ be an $x_0$-minimal submonoid of $\,\Self{X}$. Then the following are equivalent:
\begin{evlist}{8pt}{6pt}
\item[(1)]
There exists a unique monoid operation $\bullet$ with $M = M_\bullet$.

\item[(2)]
The mapping $\Phi_M$ is injective (and thus by Lemma~\ref{lemma_tm_21} bijective). 

\item[(3)]
The mapping $\Phi_{Z_M}$ is surjective (and thus by Lemma~\ref{lemma_tm_31} bijective). 

\item[(4)]
There exists a generator $A$ of $M$ and a subset $A'$ of $Z_M$  with $\Phi(A') = \Phi(A)$.

\item[(5)]
There exists a bijective mapping $\theta : M \to Z_M$ with $\Phi_{Z_M} \circ \theta = \Phi_M$.
\end{evlist}

Moreover, if one (and thus all) of these statements holds and $\bullet$ is the unique monoid operation with 
$M = M_\bullet$ then $Z_M = M_{\bullet'}$.
Also, if $A$ is a generator of $M$ and $A'$ is a subset of $Z_M$  with $\Phi(A') = \Phi(A)$ then
$A'$ is a generator of $Z_M$.
\end{theorem}

\proof This is broken up into various parts below.
\eop

Note that if $M$ is a commutative $x_0$-minimal submonoid of $\,\Self{X}$ -- and thus by Lemma~\ref{lemma_tm_21} 
$\Phi_M$ is surjective -- then $\Phi_{Z_M}$ is also surjective, since $N \subset Z_N$ holds whenever $N$ is a 
commutative submonoid. Hence by Theorem~\ref{theorem_tm_21} there exists a unique monoid operation $\bullet$ such 
that $M = M_\bullet$. Moreover, $M_\bullet \subset Z_{M_\bullet} = M_{\bullet'}$ and thus $M_\bullet = M_{\bullet'}$,
since by Theorem~\ref{theorem_tm_cayley}~(3) $\Phi_{M_\bullet}$ and $\Phi_{M_{\bullet'}}$ are both bijections.
Therefore by Lemma~\ref{lemma_tm_11} $\bullet = \bullet'$, which shows that $(X,\bullet,x_0)$ is commutative.

\begin{proposition}\label{prop_tm_11}
Let $M$ be a submonoid of $\,\Self{X}$. Then a monoid operation $\bullet$ on $(X,x_0)$ with
$M = M_\bullet$ exists if and only if $\Phi_M$ is a bijection.
Moreover, if $\bullet$ exists then it is unique.
\end{proposition}

\proof 
If $M = M_\bullet$ for some monoid operation $\bullet$ then by Theorem~\ref{theorem_tm_cayley}~(3) the mapping 
$\Phi_M$ is a bijection. Suppose conversely that $\Phi_M : M \to X$ is a bijection. 
There then exists a unique binary relation $\bullet$ on $X$ such that
\[ \Phi_M(u_1) \bullet \Phi_M(u_2) = \Phi_M(u_1 \circ u_2)\]
for all $u_1,\,u_2 \in M$. The operation $\bullet$ is associative since $\circ$ is: If 
$x_1,\,x_2,\,x_3 \in X$ and $u_1,\,u_2,\,u_3 \in M$ are such that $x_j = \Phi_M(u_j)$ 
for $j = 1,\,2,\,3$ then
\begin{eqnarray*}
(x_1 \bullet x_2) \bullet x_3 &=& (\Phi_M(u_1) \bullet \Phi_M(u_2)) \bullet \Phi_M(u_3) 
= \Phi_M(u_1 \circ u_2) \bullet \Phi_M(u_3)\\
&=& \Phi_M( (u_1 \circ u_2) \circ u_3) = \Phi_M( u_1 \circ (u_2 \circ u_3))\\
&=& \Phi_M(u_1) \bullet \Phi_M(u_2 \circ u_3)
= \Phi_M(u_1) \bullet (\Phi_M(u_2) \bullet \Phi_M(u_3))\\ 
&=& x_1 \bullet (x_2 \bullet x_3) \;.
\end{eqnarray*}
Also, $x_0$ is the unit  for $\bullet$, since if $x \in X$ and $u \in M$ is such that $x = \Phi_M(u)$ then
$x \bullet x_0 = \Phi_M(u) \bullet \Phi_M(\id_X) = \Phi_M(u \circ \id_X) = \Phi_M(u) =  x$, 
and in the same way $x_0 \bullet x = x$. Therefore $(X,\bullet,x_0)$ is a monoid. 
We next show that $M$ is exactly the set of translations in 
$(X,\bullet,x_0)$. Let $u \in M$ and $x \in X$. Since $\Phi_M$ is surjective there exists $v \in M$ with 
$x = \Phi_M(v) = v(x_0)$ and thus 
\begin{eqnarray*}
 u(x) = u(\Phi_M(v)) &=& u(v(x_0)) = (u \circ v)(x_0)\\
 &=& \Phi_M(u \circ v) = \Phi_M(u) \bullet \Phi_M(v) = \Phi_M(u) \bullet x = u(x_0) \bullet x
\end{eqnarray*}
which shows $u$ is a translation in $(X,\bullet,x_0)$. Suppose conversely that $u \in \Self{X}$ is a translation
in $(X,\bullet,x_0)$; again since $\Phi_M$ is surjective there exists $v \in M$ with 
$v(x_0) = \Phi_M(v) = u(x_0)$ and then $u(x) = u(x_0) \bullet x = v(x_0) \bullet x = v(x)$ for all 
$x \in X$ (since $v$ is also a translation in $(X,f,x_0)$). Thus $u = v \in M$. Therefore $M$ consists of 
exactly the translations in $(X,\bullet,x_0)$, and so by Theorem~\ref{theorem_tm_cayley}~(2)
$M = M_\bullet$.
The final statement (concerning the uniqueness of $\bullet$) follows immediately from
Lemma~\ref{lemma_tm_11}.
\eop

\begin{lemma}\label{lemma_tm_41}
For each monoid operation $\bullet$ on $(X,x_0)$ the centraliser of $M_\bullet$ in $\,\Self{X}$ is the 
submonoid $M_{\bullet'}$,i.e., $Z_{M_\bullet} = M_{\bullet'}$. 
\end{lemma}

\proof 
Consider $u \in M_\bullet$ and $v \in M_{\bullet'}$; then by Theorem~\ref{theorem_tm_cayley}~(2) $u$ is a 
translation in $(X,\bullet,x_0)$ and $v$ a translation in $(X,\bullet',x_0)$, hence for all $x \in X$
\begin{eqnarray*}
(u\circ v)(x) &=& u(v(x)) = u(x_0) \bullet v(x)\\ 
&=& u(x_0) \bullet (v(x_0) \bullet' x) = u(x_0) \bullet (x \bullet v(x_0)) 
= (u(x_0) \bullet x) \bullet v(x_0) \\
&=& u(x) \bullet v(x_0) = v(x_0) \bullet' u(x) = v(u(x)) = (v \circ u)(x) 
\end{eqnarray*}
and so $u \circ v = v \circ u$. Therefore $v \in Z_{M_\bullet}$, which implies $M_{\bullet'} \subset Z_{M_\bullet}$.
Now consider $u \in Z_{M_\bullet}$, and so $u \circ v = v \circ u$ for all $v \in M_\bullet$; we show that
$u(x) = u(x_0) \bullet' x$ for all $x \in X$, which will imply that 
$u = \Psi_{\bullet'}(u(x_0)) \in M_{\bullet'}$. Thus let $x \in X$; then $x = v(x_0)$ for some $v \in M_\bullet$, 
since $\Phi_{M_\bullet}$ is surjective, and hence
\[ u(x) = u(v(x_0)) = v(u(x_0)) = v(x_0) \bullet u(x_0) = x \bullet u(x_0) = u(x_0) \bullet' x\;.\]
Hence $Z_{M_\bullet} \subset M_{\bullet'}$.
\eop

\begin{proposition}\label{prop_tm_21}
Let $M$ be an $x_0$-minimal submonoid of $\,\Self{X}$. Then the mapping $\Phi_M$ is injective (and thus by
Lemma~\ref{lemma_tm_21} bijective) if and only if $\Phi_{Z_M}$ is surjective (and thus by Lemma~\ref{lemma_tm_31} 
bijective). Moreover, in this case $Z_M = M_{\bullet'}$, where $\bullet$ is the unique monoid operation such that 
$M = M_\bullet$ (given by Proposition~\ref{prop_tm_11}).
\end{proposition}

\proof 
Suppose first $\Phi_{Z_M}$ is surjective.
Let $u_1,\,u_2 \in M$ with $\Phi_M(u_1) = \Phi_M(u_2)$, i.e., with $u_1(x_0) = u_2(x_0)$. 
Then the set $X_0 = \{ x \in X : u_1(x) = u_2(x) \}$ contains $x_0$, and it is $Z_M$-invariant, since if 
$u_1(x) = u_2(x)$ then for all $v \in Z_M$
\[u_1(v(x)) = v(u_1(x)) = v(u_2(x)) = u_2(v(x))\;.\]
Hence $X_0 = X$, since by Lemma~\ref{lemma_tm_21} $Z_M$ is $x_0$-minimal, i.e., $u_1 = u_2$, which implies 
that $\Phi_M$ is injective. (Note that this the same as the proof of Lemma~\ref{lemma_tm_31}, but with the roles
of $M$ and $Z_M$ reversed).
Suppose now that $\Phi_M$ is injective and thus bijective. Then by Proposition~\ref{prop_tm_11}
there exists a monoid operation $\bullet$ with $M = M_\bullet$ and so by Lemma~\ref{lemma_tm_41}
$Z_M = M_{\bullet'}$. Therefore by Proposition~\ref{prop_tm_11} $\Phi_{Z_M}$ is bijective.
\eop

\begin{proposition}\label{prop_tm_31}
Let $M$ be an $x_0$-minimal submonoid of $\,\Self{X}$ and let $A$ be a generator of $M$. Suppose 
there exists a subset $A'$ of $Z_M$ with $\Phi(A) = \Phi(A')$. Then $\Phi_{Z_M}$ is surjective (and thus by
Lemma~\ref{lemma_tm_31} bijective) and $A'$ is a generator of $Z_M$. 
\end{proposition}

\proof 
Let $N$ be any submonoid of $Z_M$ containing $A'$ and put $X_0 = \Phi(N)$. Then the set $X_0$ is $A$-invariant:
Let $x = v(x_0) \in X_0$ (with $v \in N$) and $u \in A$, and so
$u \circ v = v \circ u$ (since $u \in M$ and $v \in Z_M$); moreover, since $\Phi(A) = \Phi(A')$
there exists $u' \in N$ with $\Phi(u') = \Phi(u)$, i.e., with $u'(x_0) = u(x_0)$. Hence
\[ u(x) = u(v(x_0)) = v(u(x_0)) = v(u'(x_0)) = (v \circ u')(x_0) = \Phi(v \circ u') \in X_0\]
and so $u(x) \in X_0$. Thus $X_0$ is $M$-invariant, since $\{ u \in \Self{X} : \mbox{$X_0$ is $u$-invariant}\}$
is a submonoid of $\Self{X}$. Moreover, $x_0 \in X_0$, since $\id_X \in N$ and $\id_X(x_0) = x_0$. Therefore 
$X_0 = X$, since $M$ is $x_0$-minimal, which shows that $\Phi_N$ is surjective.
But $N \subset Z_M$ and so $\Phi_{Z_M}$ is surjective. Therefore by Lemma~\ref{lemma_tm_31}
$\Phi_{Z_M}$ is bijective, which is only possible if $N = Z_M$. Moreover, taking $N$ to be the least 
submonoid containing $A'$ implies that $A'$ is a generator of $Z_M$.
\eop

\textit{Proof of Theorem~\ref{theorem_tm_21}:}\enskip
(1) $\Leftrightarrow$ (2) is Proposition~\ref{prop_tm_11}.

(2) $\Leftrightarrow$ (3):\enskip This is a part of Proposition~\ref{prop_tm_21}.

(4) $\Rightarrow$ (3):\enskip This follows from Proposition~\ref{prop_tm_31}.

(5) $\Rightarrow$ (4):\enskip This is clear, since a generator $A$ of $M$ exists (for example, $M$ itself),
and then $\Phi(A) = \Phi(A')$ with $A' = \theta(A)$.

(1) $\Rightarrow$ (5):\enskip 
Let $\bullet$ be the unique monoid operation with $M = M_\bullet$. Then by Lemma~\ref{lemma_tm_41}
$Z_M = M_{\bullet'}$ and thus by Theorem~\ref{theorem_tm_cayley}
$\theta = \Psi_{\bullet'} \circ \Phi_M : M \to Z_M$ is a bijection.

Finally, suppose that one (and thus all) of the statements holds and let $\bullet$ be the unique monoid operation 
with $M = M_\bullet$. Then by Lemma~\ref{lemma_tm_41} $Z_M = M_{\bullet'}$. Moreover,
if $A$ is a generator of $M$ and $A'$ is a subset of $Z_M$  with $\Phi(A') = \Phi(A)$ then by 
Proposition~\ref{prop_tm_31} $A'$ is a generator of $Z_M$.

This completes the proof of Theorem~\ref{theorem_tm_21}.
\eop

We now look at how Theorem~\ref{theorem_tm_21} can be applied to the situation considered in Section~\ref{lsalg}, 
and first note the following:

\begin{lemma}\label{lemma_tm_51}
Let $\Lambda = (X,f,x_0)$ be a minimal $\mathcal{L}_S$-algebra and let $\bullet$ be a monoid operation on $(X,x_0)$.
Then $\bullet$ is $\Lambda$-compatible if and only if $M_f = M_\bullet$.
\end{lemma}

\proof 
Since $\Lambda$ is minimal the submonoid $M_f$ is $x_0$-minimal.
Suppose first that $\bullet$ is $\Lambda$-compatible. Then by Theorem~\ref{theorem_tm_cayley}~(2)
$f_s \in M_\bullet$ for each $s \in S$ and so $M_f \subset M_\bullet$.
But by Lemma~\ref{lemma_tm_21} the mapping $\Phi_{M_f}$ is surjective and by 
Theorem~\ref{theorem_tm_cayley}~(3) the mapping $\Phi_{M_\bullet}$ is bijective, which implies that
$M_f = M_\bullet$. Conversely, if $M_f = M_\bullet$ then by
Theorem~\ref{theorem_tm_cayley}~(2) $\bullet$ is $\Lambda$-compatible. 
\eop

Let $\Lambda = (X,f,x_0)$ be a minimal $\mathcal{L}_S$-algebra and suppose there exists a reflection $f'$ of $f$ in 
$\Lambda$. We apply Theorem~\ref{theorem_tm_21} to show that there exists a $\Lambda$-compatible monoid operation 
$\bullet$. As before put $f_S = \{ u \in \Self{X} : \mbox{$u = f_s$ for some $s \in S$} \}$, and so $f_S$ is a 
generator of $M_f$, and put $f'_S = \{ v \in \Self{X} : \mbox{$v = f'_s$ for some $s \in S$} \}$. Then 
$\Phi(f_S) = \Phi(f'_S)$, since $f'_s(x_0) = f_s(x_0)$ for each $s \in S$. Moreover, $f'_S \subset Z_{f_S}$, since 
$f_s \circ f'_t = f'_f \circ f_s$ for all $s,\,t \in S$, and it follows that $f'_S \subset Z_{M_f}$ (since if $A$ is 
a generator of a monoid $M$ then it is easy to see that $Z_A = Z_M$). Therefore by Theorem~\ref{theorem_tm_21} 
((4) $\Rightarrow$ (1)) there exists a unique monoid operation $\bullet$ on $(X,x_0)$ with $M_f = M_\bullet$, and 
so by Lemma~\ref{lemma_tm_51} $\bullet$ is $\Lambda$-compatible.

It also follows from Theorem~\ref{theorem_tm_21} that $Z_{M_f} = M_{\bullet'}$ and that $f'_S$ is a generator of 
$M_{\bullet'}$, i.e., $M_{f'} = M_{\bullet'}$. Hence by Theorem~\ref{theorem_tm_cayley}~(3) $\Phi_{M_{f'}}$ is
bijective which, together with Lemma~\ref{lemma_tm_21}, implies that $\Lambda' = (X,f',x_0)$ is minimal.
(This was the final, and only non-trivial, statement in Proposition~\ref{prop_lsalg_21}.)

We next consider a result which corresponds to Proposition~\ref{prop_lsalg_11}.

\begin{proposition}\label{prop_tm_41} 
Let $\bullet$ be a monoid operation on $(X,x_0)$. Then:
\begin{evlist}{8pt}{6pt}
\item[(1)]
The monoid $(X,\bullet,x_0)$ is  commutative if and only if $M_\bullet$ is commutative.

\item[(2)]
The monoid $(X,\bullet,x_0)$ obeys the left cancellation law if and only each mapping
\phantom{xxx}in $M_\bullet$ is injective.

\item[(3)]
The monoid $(X,\bullet,x_0)$ is a group if and only if each mapping in $M_\bullet$ is 
\phantom{xxx}surjective, which is the case if and only if each mapping in $M_\bullet$ is a bijection.
\end{evlist}

\end{proposition}

\proof
These  all use the fact that the monoids $(X,\bullet,x_0)$ and $(M_\bullet,\circ,\id_X)$ are isomorphic
(which was established in Theorem~\ref{theorem_tm_cayley}). In particular, the monoid $(X,\bullet,x_0)$ is 
commutative if and only $M_\bullet$ is commutative, which is (1).
For parts (2) and (3) we need the following fact:

\begin{lemma}\label{lemma_tm_61}
Let $M$ be a submonoid of $\,\Self{X}$ for which the mapping $\Phi_M$ is bijective. Then:
\begin{evlist}{8pt}{6pt}
\item[(1)]
The monoid $(M,\circ,\id_X)$ obeys the left cancellation law if and only if each 
\phantom{xxx}mapping in $M$ is injective. 

\item[(2)]
The monoid $(M,\circ,\id_X)$ is a group  if and only if each mapping in $M$ is 
\phantom{xxx}surjective, which is the case if and only if each mapping in $M$ is a
bijection. 
\end{evlist}
\end{lemma}

\proof
(1)\enskip
Suppose $(M,\circ,\id_X)$ obeys the left cancellation law. Let $u \in M$ and $x_1,\,x_2 \in X$ with 
$u(x_1) = u(x_2)$. Then there exist $u_1,\,u_2 \in M$ with $\Phi_M(u_1) = x_1$ and $\Phi_M(u_2) = x_2$ 
(since $\Phi_M$ is surjective), and hence 
\begin{eqnarray*}
\Phi_M(u\circ u_1) &=& (u \circ u_1)(x_0) = u(u_1(x_0)) = u(\Phi_M(u_1)) = u(x_1)\\ 
&=& u(x_2) = u(\Phi_M(u_2)) = u(u_2(x_0)) = (u \circ u_2)(x_0) = \Phi_M(u\circ u_2)\;.
\end{eqnarray*} 
It follows that $u \circ u_1 = u \circ u_2$ (since $\Phi_M$ is injective) and so $u_1 = u_2$. In particular 
$x_1 = x_2$, which shows that $u$ is injective. The converse is immediate, since if $u \in M$ is injective and 
$u \circ u_1 = u \circ u_2$ then $u_1 = u_2$.

(2)\enskip
Suppose that each mapping in $M$ is surjective. Let $u \in M$; there then exists $x \in X$ such that $u(x) = x_0$ 
(since $u$ is surjective) and there exists $v \in M$ with $\Phi_M(v) = x$ (since $\Phi_M$ is surjective). 
It follows that
\[  \Phi_M(u \circ v) = (u \circ v)(x_0) = u(v(x_0)) = u(\Phi_M(v)) = u(x) = x_0 = \Phi_M(\id_X) \]
and therefore $u \circ v = \id_X$, since $\Phi_M$ is injective. For each $u \in M$ there thus exists $v \in M$ 
with $u \circ v = \id_X$, and therefore (as in the proof of Proposition~\ref{prop_lsalg_11}~(3)) $(M,\circ,\id_X)$ 
is a group. Suppose conversely  $(M,\circ,\id_X)$ is a group. Then for each $u \in M$ there exists
$v \in M$ with $u \circ v = v \circ u = \id_X$ and so $u$ is a bijection.
\eop

We apply Lemma~\ref{lemma_tm_61} to obtain parts (2) and (3) of Proposition~\ref{prop_tm_41}:

(2):\enskip
The monoid $(X,\bullet,x_0)$ obeys the left cancellation law if and only if $(M,\circ,\id_X)$ does, which by 
Lemma~\ref{lemma_tm_61}~(1) is the case if and only if if each $v \in M$ is injective.

(3):\enskip
The monoid $(X,\bullet,x_0)$ is a group if and only if $(M,\circ,\id_X)$ is, which by 
Lemma~\ref{lemma_tm_61}~(2) is the case if and only if if each mapping in $M$ is surjective,
and this is the case  if if each mapping in $M$ is a bijection.

This completes the proof of Proposition~\ref{prop_tm_41}.
\eop

The statements in Proposition~\ref{prop_lsalg_11} can all be deduced from
Proposition~\ref{prop_tm_41}. Let $\Lambda = (X,f,x_0)$ be a minimal $\mathcal{L}_S$-algebra 
for which there exists a $\Lambda$-compatible monoid operation $\bullet$, and so by Lemma~\ref{lemma_tm_61}
$M_f = M_\bullet$. By Proposition~\ref{prop_tm_41}~(1) $(X,\bullet,x_0)$ is commutative
if and only if $M_f$ is, and by the lemma following this is the case if and only if $\Lambda$ is commutative.

\begin{lemma}\label{lemma_tm_71}
If a submonoid $M$ of $\,\Self{X}$ has a commutative generator $A$ (meaning that $u \circ v = v \circ u$ for 
all $u,\,v \in A$) then $M$ is  commutative.
\end{lemma}

\proof 
For each $u \in \Self{X}$ the set $C_u = \{\,v \in \Self{X} : \mbox{$v \circ u = u \circ v$}\, \}$ is a 
submonoid, since $\id_X \circ u = u = u \circ \id_X$ and if $u_1,\,u_2 \in C_u$ then 
\[ (u_1 \circ u_2) \circ u = u_1 \circ u_2 \circ u = u_1 \circ u \circ u_2 
= u \circ u_1 \circ u_2 = u \circ (u_1 \circ u_2)\;. \] 
Now $A \subset C_u$ for all $u \in A$, since $A$ is a commutative subset, and thus $M \subset C_u$ for all 
$u \in A$, i.e., $v \circ u = u \circ v$ for all $u \in A$, $v \in M$. But this also says that $A \subset C_v$ 
for each $v \in M$, which implies that $M \subset C_v$ for each $v \in M$, and hence shows that
$v \circ u = u \circ v$ for all $u,\, v \in M$. \eop

Similarly, by Proposition~\ref{prop_tm_41}~(2) $(X,\bullet,x_0)$ obeys the left cancellation law if and only
if each mapping in $M_f$ is injective, which is the case if and only if $f_s$ is injective for each $s \in S$,
since $\{ u \in \Self{X} : \mbox{$u$ is injective} \}$ is a submonoid of $\Self{X}$.
In the same way Proposition~\ref{prop_tm_41}~(3) implies $(X,\bullet,x_0)$ is a group if and only if $f_s$ is 
surjective for each $s \in S$, and which is the case if and only if $f_s$ is bijective for each $s \in S$, since 
$\{ u \in \Self{X} : \mbox{$u$ is surjective} \}$ is also a submonoid of $\Self{X}$.


%% file: initial.tex
\startsection{Initial $\mathcal{L}_S$-algebras}
\label{initial}

In this section we look at in more detail at the $\mathcal{L}_S$-algebra
$\Lambda = (S^*,\triangleleft,\varepsilon)$ of `real' lists of elements from $S$. Recall that 
$S^* = \bigcup_{n \ge 0} S^n$, with the element $(s_1,\ldots,s_n)$ of $S^n$ usually written as 
$s_1\ \cdots\ s_n$, $\varepsilon$ is the single element in $S^0$ and 
$\triangleleft_s(s_1\ \cdots\ s_n) = s\ s_1\ \cdots\ s_n$, with 
$\triangleleft_s(\varepsilon) = s \in S^1 \subset S^*$. This $\mathcal{L}_S$-algebra is minimal and the 
concatenation operation $\bullet$ given by
\[s_1\ \cdots\ s_m \,\bullet\, s'_1\ \cdots\ s'_n = s_1\ \cdots\ s_m \ s'_1\ \cdots\ s'_n\]
is the unique $\Lambda$-compatible monoid operation. Moreover, the reflection $\triangleleft'$ of 
$\triangleleft$ in $\Lambda$ is the mapping such that $\triangleleft'_s$ is the operation of adding the 
element $s$ to the end of a list, i.e., $\triangleleft'_s(s_1\ \cdots\ s_n) = s_1\ \cdots\ s_n\ s$. Another 
mapping which plays a role here is the mapping $r : S^* \to S^*$ which reverses a list, and so
\[r(s_1\ \cdots\ s_m) = s_m\ \cdots\ s_1\]
for each list $s_1\ \cdots\ s_m$. Thus $r(\varepsilon) = \varepsilon$ and both
$\triangleleft'_s \circ r = r \circ \triangleleft_s$ and $\triangleleft_s \circ r = r \circ \triangleleft'_s$ 
hold for each $s \in S$. Moreover, $r \circ r = \id_{S^*}$, i.e., reversing a list twice ends up with the original 
list.

Now it might appear that there is not much more to say about this $\mathcal{L}_S$-algebra with respect to the 
topics we have been considering. However, there are a couple of points which are not very satisfactory. The 
first concerns the implicit use of properties of the natural numbers in defining 
$(S^*,\triangleleft,\varepsilon)$. For example, the definition of the set $S^n$ involves the segment 
$\{1,2,\ldots,n\}$, whose properties are usually taken for granted, but which are not so trivial to establish 
starting with the Peano axioms. We would prefer to avoid this dependence, in particular since, except for their 
appearance in some of the examples, the natural numbers have played no role in these notes.

The second point is that some explanation is needed for why $(S^*,\triangleleft,\varepsilon)$ behaves like it 
does, and the reason is that $(S^*,\triangleleft,\varepsilon)$ is an initial $\mathcal{L}_S$-algebra. This 
basic fact is well-known and is usually taught in some form in most introductory computer science courses.
We present this topic here, tying it in with the results from the previous sections and without in any way 
making use of the natural numbers.

We start by introducing the structure preserving mappings between $\mathcal{L}_S$-algebras. If $(X,f,x_0)$ and 
$(Y,g,y_0)$ are $\mathcal{L}_S$-algebras then a mapping $\pi : X \to Y$ is called a 
\definition{morphism from $(X,f,x_0)$ to $(Y,g,y_0)$} if $\pi(x_0) = y_0$ and $g_s \circ \pi = \pi \circ f_s$ 
for all $s \in S$. This will also be indicated by stating that $\pi : (X,f,x_0) \to (Y,g,y_0)$ is a morphism.

\begin{lemma}\label{lemma_initial_w1}
(1)\enskip
For each $\mathcal{L}_S$-algebra $(X,f,x_0)$ the identity mapping $\id_X$ is a morphism from $(X,f,x_0)$ to 
$(X,f,x_0)$. 

(2)\enskip
If $\pi : (X,f,x_0) \to (Y,g,y_0)$ and $\sigma : (Y,g,y_0) \to (Z,h,z_0)$ are morphisms then $\sigma \circ \pi$ 
is a morphism from $(X,f,x_0)$ to $(Z,h,z_0)$. 
\end{lemma}

\proof 
(1)\enskip
This is clear, since $\id_X(x_0) = x_0$ and $f_s \circ \id_X = f_s = \id_X \circ f_s$ for all $s \in S$.

(2)\enskip
This follows since $(\sigma \circ \pi)(x_0) = \sigma(\pi(x_0)) = \sigma(y_0) = z_0$ and
\[ h_s\circ (\sigma \circ \pi) =  (h_s\circ \sigma) \circ \pi =  (\sigma \circ g_s) \circ \pi 
=  \sigma \circ (g_s \circ \pi) =  \sigma \circ (\pi \circ f_s) = (\sigma \circ \pi) \circ f_s 
\]
for all $s \in S$. \eop

If $\pi : (X,f,x_0) \to (Y,g,y_0)$ is a morphism then clearly $\pi \circ \id_X = \pi = \id_Y \circ \pi$, and 
if $\pi,\,\sigma$ and $\tau$ are morphisms for which the compositions are defined then
$(\tau \circ \sigma) \circ \pi = \tau \circ (\sigma \circ \pi)$. This means that $\mathcal{L}_S$-algebras are 
the objects of a concrete category, whose morphisms are those defined above.

Note that if $(X,f,x_0)$ is minimal then for each $\mathcal{L}_S$-algebra $(Y,g,y_0)$ there can be at most one 
morphism $\pi : (X,f,x_0) \to (Y,g,y_0)$, since if $\pi_1$ and $\pi_2$ are two such morphisms then the set 
$X_0 = \{ x \in X : \pi_1(x) = \pi_2(x) \}$ contains $x_0$ and is easily seen to be $f$-invariant. Thus 
$X_0 = X$, since $(X,f,x_0)$ is minimal, i.e., $\pi_1 = \pi_2$.

An \definition{isomorphism} is a morphism $\pi : (X,f,x_0) \to (Y,g,y_0)$ for which there exists a morphism 
$\sigma : (Y,g,y_0) \to (X,f,x_0)$ such that $\sigma \circ \pi = \id_X$ and $\pi \circ \sigma = \id_Y$. In this 
case $\sigma$ is uniquely determined by $\pi$: If $\sigma' : (Y,g,y_0) \to (X,f,x_0)$ is also a morphism with 
$\sigma' \circ \pi = \id_X$ and $\pi \circ \sigma' = \id_Y$ then
\[   \sigma' = \sigma' \circ \id_Y = \sigma' \circ (\pi \circ \sigma) 
= (\sigma' \circ \pi) \circ \sigma = \id_X \circ \sigma = \sigma\;. \]
The morphism $\sigma$ is called the \definition{inverse} of $\pi$.

\begin{lemma}\label{lemma_initial_w2}
A morphism $\pi : (X,f,x_0) \to (Y,g,y_0)$ is an isomorphism if and only if the mapping $\pi : X \to Y$ is a 
bijection; in this case the inverse morphism is the inverse mapping $\pi^{-1} : Y \to X$.
\end{lemma}

\proof 
If $\sigma \circ \pi = \id_X$ and $\pi \circ \sigma = \id_Y$ then $\pi$ is a bijection and $\sigma$ is the 
inverse mapping $\pi^{-1} : Y \to X$. It thus remains to show that if $\pi$ is a bijection then the inverse 
mapping $\pi^{-1} : Y \to X$ defines a morphism from $(Y,g,y_0)$ to $(X,f,x_0)$. Let $y \in Y$; then there 
exists a unique $x \in X$ with $y = \pi(x)$ and thus
\[ f_s(\pi^{-1}(y))) = f_s(x) = \pi^{-1}(\pi(f_s(x))) = \pi^{-1}(g_s(\pi(x)))= \pi^{-1}(g_s(y))\;,
\]
and this implies that $f_s \circ \pi^{-1} = \pi^{-1} \circ g_s$ for all $s \in S$. Moreover $\pi^{-1}(y_0) = x_0$, 
since $\pi(x_0) = y_0$, and therefore $\pi^{-1} : (Y,g,y_0) \to (X,f,x_0)$ is a morphism.
\eop

The $\mathcal{L}_S$-algebras $(X,f,x_0)$ and $(Y,g,y_0)$ are said to be \definition{isomorphic} if there exists 
an isomorphism $\pi : (X,f,x_0) \to (Y,g,y_0)$. Being isomorphic clearly defines an equivalence relation on the 
class of all $\mathcal{L}_S$-algebras.

An $\mathcal{L}_S$-algebra $(X,f,x_0)$ is said to be \definition{initial} if for each $\mathcal{L}_S$-algebra 
$(Y,g,y_0)$ there exists a unique morphism $\pi : (X,f,x_0) \to (Y,g,y_0)$.

The following simple fact about initial objects holds in any category:

\begin{lemma}\label{lemma_initial_w3}
If $(X,f,x_0)$ and $(Y,g,y_0)$ are initial $\mathcal{L}_S$-algebras then the unique morphism
$\pi : (X,f,x_0) \to (Y,g,y_0)$ is an isomorphism. In particular, $(X,f,x_0)$ and $(Y,g,y_0)$ are isomorphic.
\end{lemma}

\proof 
Since $(Y,g,y_0)$ is initial there exists a unique morphism $\sigma$ from $(Y,g,y_0)$ to $(X,f,x_0)$ and then 
by Lemma~\ref{lemma_initial_w1}~(2) $\sigma \circ \pi$ is a morphism from $(X,f,x_0)$ to $(X,f,x_0)$. But 
$(X,f,x_0)$ is initial and so there is a unique such morphism, which by Lemma~\ref{lemma_initial_w1}~(1) is 
$\id_X$, and hence $\sigma \circ \pi = \id_X$. In the same way (reversing the roles of $(X,f,x_0)$ and 
$(Y,g,y_0)$) it follows that $\pi \circ \sigma = \id_Y$ and therefore $\pi$ is an isomorphism.
\eop

An $\mathcal{L}_S$-algebra $(X,f,x_0)$ will be called \definition{unambiguous} if the mapping $f_s$ is 
injective for each $s \in S$ and the sets $f_s(X)$, $s \in S$, are disjoint and 
$x_0 \notin \bigcup_{s\in S} f_s(X)$.

Note that the $\mathcal{L}_S$-algebra $(S^*,\triangleleft,\varepsilon)$ is both minimal and unambiguous.

\begin{theorem}\label{theorem_initial_w1}
There exists an initial $\mathcal{L}_S$-algebra, and an $\mathcal{L}_S$-algebra is initial if and only if it 
is minimal and unambiguous.
\end{theorem}

The second statement in Theorem~\ref{theorem_initial_w1} is often expressed by computer scientists by saying 
that the initial objects are characterised as having  \textit{no junk} (being minimal) and
\textit{no confusion} (being unambiguous).

\begin{theorem}\label{theorem_initial_w2}
Let $\Lambda = (X,f,x_0)$ be an initial $\mathcal{L}_S$-algebra. Then:
\begin{evlist}{8pt}{6pt}
\item[(1)]
There exists a unique $\Lambda$-compatible monoid operation $\,\bullet$, and the monoid
\phantom{xxx}$(X,\bullet,x_0)$ obeys both the left and right cancellation laws.

\item[(2)]
If $f'$ is the reflection of $f$ in $\Lambda$ 
then the $\mathcal{L}_S$-algebra $\Lambda' = (X,f',x_0)$ is also 
\phantom{xxx}initial.

\item[(3)]
If $r : (X,f,x_0) \to (X,f',x_0)$ is the unique morphism (so by Lemma~\ref{lemma_initial_w3} $r$ 
\phantom{xxx}is an isomorphism) then $r$ is also the unique isomorphism from $(X,f',x_0)$ 
\phantom{xxx}to $(X,f,x_0)$ and $r \circ r = \id_X$.
\end{evlist}
\end{theorem}

We now start preparing for the proofs of Theorems \ref{theorem_initial_w1} and \ref{theorem_initial_w2}.

\begin{lemma}\label{lemma_initial_w4}
Let $(X,f,x_0)$ be a minimal $\mathcal{L}_S$-algebra. Then for each $x \in X \setminus \{x_0\}$ there exists 
$x' \in X$ and $s \in S$ so that $x = f_s(x')$.
\end{lemma}

\proof 
Let $X_0$ be the subset of $X$ consisting of $x_0$ together with all elements of the form $f_s(x)$ with 
$s \in S$ and $x \in X$. Then $X_0$ is clearly $f$-invariant and it contains $x_0$ and hence $X_0 = X$, since 
$(X,f,x_0)$ is minimal. \eop

Lemma~\ref{lemma_initial_w4} shows that if $(X,f,x_0)$ is a minimal unambiguous $\mathcal{L}_S$-algebra then 
for each element $x \in X \setminus \{x_0\}$ there exists a unique $s \in S$ and a unique $x' \in X$ such that 
$x = f_s(x')$.

\begin{lemma}\label{lemma_initial_w5}
Let $(X,f,x_0)$ be any $\mathcal{L}_S$-algebra, let $X^0$ be the least $f$-invariant subset of $X$ containing 
$x_0$ and for each $s \in S$ let $f^0_s$ be the restriction of $f_s$ to $X^0$, considered as a mapping from 
$X^0$ to itself. Then the $\mathcal{L}_S$-algebra $(X^0,f^0,x_0)$ is minimal.
\end{lemma}

\proof 
An $f^0$-invariant subset $X'$ of $X^0$ containing $x_0$ is  also an $f$-invariant subset of $X$ containing 
$x_0$ and so $X^0 \subset X'$. Thus $X' = X^0$, which implies that the only $f^0$-invariant subset of $X^0$ 
containing $x_0$ is $X^0$ itself. Therefore $(X^0,f^0,x_0)$ is a minimal $\mathcal{L}_S$-algebra.
\eop

\begin{lemma}\label{lemma_initial_w6}
Let $(X,f,x_0)$ be a minimal $\mathcal{L}_S$-algebra, $(Y,g,y_0)$ an unambiguous $\mathcal{L}_S$-algebra and 
suppose there exists a morphism $\pi : (X,f,x_0) \to (Y,g,y_0)$. Then $\pi$ is injective and $(X,f,x_0)$ is 
unambiguous.
\end{lemma}

\proof 
We first show that $\pi$ is injective, and so consider the set
\[  X_0 = \{ x \in X : \mbox{$x' = x$ whenever $x' \in X$ with $\pi(x') = \pi(x)$}  \}\;. \]
If $x' \in X \setminus \{x_0\}$ then by Lemma~\ref{lemma_initial_w4} there exists $s \in S$ and $x'' \in X$ 
with $x' = f_s(x'')$ and so $\pi(x') = \pi(f_s(x'')) = g_s(\pi(x'')) \ne y_0 = \pi(x_0)$. Hence $x_0 \in X_0$.
Also $X_0$ is $f$-invariant: Let $x \in X_0$ and $s \in S$ and suppose $\pi(f_s(x)) = \pi(x')$ for some 
$x' \in X$. Then $x' \ne x_0$, since $\pi(f_s(x)) = g_s(\pi(x)) \ne y_0 = \pi(x_0)$, and thus by 
Lemma~\ref{lemma_initial_w4} there exists $t \in S$ and $x'' \in X$ with $x' = f_t(x'')$. It follows that 
$g_s(\pi(x)) = \pi(f_s(x)) = \pi(x') = \pi(f_t(x'')) = g_t(\pi(x''))$, which is only possible if $s = t$ and 
$\pi(x) = \pi(x'')$. Hence $x = x''$, since $x \in X_0$, and so $x' = f_s(x)$, which means that 
$f_s(x) \in X_0$. Therefore $X_0 = X$, since $(X,f,x_0)$ is minimal, i.e., $\pi$ is injective.

It follows immediately that $f_s$ is injective for each $s \in S$, since $\pi$ and $g_s$ are injective and 
$\pi \circ f_s = g_s \circ \pi$. Moreover $\pi(f_s(x)) = g_s(\pi(x)) \ne y_0 = \pi(x_0)$ and so 
$f_s(x) \ne x_0$ for all $x \in X$, $s \in S$. Finally, $f_s(x) = f_t(x')$ can only hold if $s = t$ and 
$x = x'$ since then $g_s(\pi(x)) = \pi(f_s(x)) = \pi(f_t(x')) = g_t(\pi(x'))$. Hence $(X,f,x_0)$ is 
unambiguous. 
\eop

\begin{lemma}\label{lemma_initial_w7}
There exists an unambiguous minimal $\mathcal{L}_S$-algebra.
\end{lemma}

\proof 
If $(X,f,x_0)$ is unambiguous and $(X^0,f^0,x_0)$ is the minimal $\mathcal{L}_S$-algebra given in 
Lemma~\ref{lemma_initial_w5} then clearly $(X^0,f^0,x_0)$ is also unambiguous. It is thus enough to show that 
an unambiguous $\mathcal{L}_S$-algebra exists.

Choose any infinite set $A$, and so there exists a proper subset $A_0$ of $A$ and a surjective mapping
$\gamma : A_0 \to A$; also let $\varepsilon$ be some element not in $S$. Now let $X$ be set of all mappings 
from $A$ to $S \cup \{\varepsilon\}$, let $x_0$ be the constant mapping with $x_0(a) = \varepsilon$ for all 
$a \in A$, and for each $s \in S$ let $f_s : X \to X$ be given by
\[   f_s(x)(a) = \left\{
            \begin{array}{cl}
                   x(\gamma(a)) &\ \mbox{if $a \in A_0$}\;,\\
                   s &\ \mbox{if $a \in A \setminus A_0$}\;.
                \end{array} \right. \]
Then the $\mathcal{L}_S$-algebra $(X,f,x_0)$ is unambiguous: $f_s(x)(a) = s \ne \varepsilon = x_0(a)$ for all 
$a \in A \setminus A_0$, and thus $x_0 \notin f_s(X)$ for each $s \in S$. In the same way, if $s \ne t$ then
$f_s(x)(a) = s \ne t = f_t(x')$ for all $a \in A \setminus A_0$, and so $f_s(X)$ and $f_t(X)$ are disjoint. 
Finally, if $f_s(x) = f_s(x')$ then $x(\gamma(a)) = x'(\gamma(a))$ for all $a \in A_0$, and hence $x = x'$, 
since $\gamma : A_0 \to A$ is surjective. This shows that $f_s$ is injective for each $s \in S$. 
\eop

Note the use of the infinite set $A$ and the mapping $\gamma : A_0 \to A$ in the above proof. If there was a 
need to be explicit we could here take $A = \Nat$, $A_0 = \Nat \setminus \{0\}$ and $\gamma : A_0 \to A$ to 
be the unique mapping with $\gamma(\textsf{s}(n)) = n$ for all $n \in A$.

\begin{lemma}\label{lemma_initial_w8}
An initial $\mathcal{L}_S$-algebra $(X,f,x_0)$ is unambiguous and minimal. 
\end{lemma}

\proof 
We first show that $(X,f,x_0)$ is minimal. Consider the minimal $\mathcal{L}_S$-algebra $(X^0,f^0,x_0)$ given 
in Lemma~\ref{lemma_initial_w5}. Then the inclusion mapping of $X^0$ in $X$ results in a morphism
$\mathrm{inc} : (X^0,f^0,x_0) \to (X,f,x_0)$ and there exists a unique morphism 
$\pi : (X,f,x_0) \to (X^0,f^0,x_0)$. Therefore by Lemma~\ref{lemma_initial_w1}~(2) there is a morphism
$\mathrm{inc}\circ \pi : (X,f,x_0) \to (X,f,x_0)$. But $\id_X$ is the unique such morphism, and so 
$\mathrm{inc} \circ \pi = \id_X$. Hence $X = \id_X(X) = \mathrm{inc}(\pi(X)) \subset X_0$, i.e., $X = X^0$,
and thus $(X,f,x_0)$ is minimal. 

It remains to show that $(X,f,x_0)$ is unambiguous. By Lemma~\ref{lemma_initial_w7} there exists an unambiguous
$\mathcal{L}_S$-algebra $(Y,g,x_0)$, so let $\pi : (X,f,x_0) \to (Y,g,y_0)$ be the unique morphism. Since 
$(X,f,x_0)$ is minimal we can apply Lemma~\ref{lemma_initial_w6}, which gives us that $(X,f,x_0)$ is unambiguous.
\eop

\begin{lemma}\label{lemma_initial_w9}
An unambiguous minimal $\mathcal{L}_S$-algebra is initial.
\end{lemma}

\proof
The proof is almost identical to one of the standard proofs of the recursion theorem. Let $(X,f,x_0)$ be an 
unambiguous minimal $\mathcal{L}_S$-algebra and $(Y,g,y_0)$ be any $\mathcal{L}_S$-algebra, and consider the 
$\mathcal{L}_S$-algebra $(X \times Y, f \times_S g, (x_0,y_0))$, where
$f \times_S g : S \times X \times Y \to X \times Y$ is given by
$(f \times_S g)(s,x,y) = (f(s,x),g(s,y))$ for all $s \in S$, $x \in X$, $y \in Y$, and so
$(f \times_S g)_s = f_s \times g_s$ for each $s \in S$. Let $Z$ be the least $(f \times_S g)$-invariant subset 
of $X \times Y$ containing $(x_0,y_0)$ and let
\[ X_0 = \{\, x \in X : \mbox{there exists exactly one $y \in Y$ such that $(x,y) \in Z$} \,\}\;. \]
It will be shown that $X_0$ is an $f$-invariant subset of $X$ containing $x_0$, which implies that $X_0 = X$,
since $(X,f,x_0)$ is minimal. We twice need the following fact: If $(x,y) \in Z \setminus \{(x_0,y_0)\}$ then 
there exists $s \in S$ and $(x',y') \in Z$ such that $(f_s(x'),g_s(y')) = (x,y)$. (This follows because
$\{(x_0,y_0)\} \cup \bigcup_{s \in S} (f_s \times g_s)(Z)$ is an $(f \times_S g)$-invariant subset of 
$X \times Y$ containing $(x_0,y_0)$ and so contains $Z$.)

The element $x_0$ is in $X_0$: Clearly $(x_0,y_0) \in Z$, so suppose also $(x_0,y) \in Z$ for some $y \ne y_0$. 
Then $(x_0,y) \in Z \setminus \{(x_0,y_0)\}$ and hence there exists $(x',y') \in Z$ and $s \in S$ with 
$(f_s(x'),g_s(y')) = (x_0,y)$. In particular $f_s(x') = x_0$, which is not possible, since $(X,f,x_0)$ is 
unambiguous. This shows that $x_0 \in X_0$.

Next let $x \in X_0$ and $s \in S$ and let $y$ be the unique element of $Y$ with $(x,y) \in Z$. Hence 
$(f_s(x),g_s(y)) = (f_s \times g_s)(x,y) \in Z$, since $Z$ is $(f \times_S g)$-invariant. Suppose also  
$(f_s(x),y') \in Z$ for some $y' \in Y$. Then $(f_s(x),y') \in Z \setminus \{(x_0,y_0)\}$, since 
$f_s(x) \ne x_0$, and so $(f_s(x),y') = (f_t(x''),g_t(y''))$ for some $t \in S$ and $(x'',y'') \in Z$. In 
particular $f_t(x'') = f_s(x)$, and this is only possible with $t = s$ and $x'' = x$. since $(X,f,x_0)$ is 
unambiguous. Therefore $y'' = y$, since $x \in X_0$, which implies $y' = g_s(y'') = g_s(y)$. This shows that 
$g_s(y)$ is the unique element $\breve{y} \in Y$ with $(f_s(x),\breve{y}) \in Z$ and in particular that 
$f_s(x) \in X_0$.

We have established that $X_0$ is an $f$-invariant subset of $X$ containing $x_0$, and so $X_0 = X$. Now define 
a mapping $\pi : X \to Y$ by letting $\pi(x)$ be the unique element of $Y$ such that $(x,\pi(x)) \in Z$ for 
each $x \in X$. Then $\pi(x_0) = y_0$, since $(x_0,y_0) \in Z$ and $\pi(f_s(x)) = g_s(\pi(x))$ for all 
$x \in X$, $s \in S$, since $(f(x),g(y)) \in Z$ whenever $(x,y) \in Z$ and so in particular
$(f_s(x),g_s(\pi(x))) \in Z$ for all $x \in X$, $s \in S$. This gives us a morphism $\pi$ from $(X,f,x_0)$ to 
$(Y,g,y_0)$, and it is easy to see that $f$ being minimal implies that $\pi$ is unique. Therefore the
$\mathcal{L}_S$-algebra $(X,f,x_0)$ is initial.
\eop

\textit{Proof of Theorem~\ref{theorem_initial_w1}:}\enskip
Lemmas \ref{lemma_initial_w7} and \ref{lemma_initial_w8} show that an $\mathcal{L}_S$-algebra is initial if 
an only if it is minimal and unambiguous and this, together with Lemma~\ref{lemma_initial_w9} also shows that 
an initial $\mathcal{L}_S$-algebra exists.
\eop

\begin{lemma}\label{lemma_initial_w10}
For each initial $\mathcal{L}_S$-algebra $\Lambda = (X,f,x_0)$ there exists a reflection $f'$ of $f$ in 
$\Lambda$.
\end{lemma}

\proof 
For each $s \in S$ there is a unique morphism $f'_s : (X,f,x_0) \to (X,f,f_s(x_0))$, and so 
$f'_s(x_0) = f_s(x_0)$ and $f_t \circ f'_s = f'_s \circ f_t$ for all $t \in S$. Thus $f'$ is a reflection of $f$ 
in $\Lambda$. 
\eop

\textit{Proof of Theorem~\ref{theorem_initial_w2}:}\enskip
Let $\Lambda = (X,f,x_0)$ be an initial $\mathcal{L}_S$-algebra; then by Lemma~\ref{lemma_initial_w8}
$\Lambda$ is minimal and by Lemma~\ref{lemma_initial_w10} there exists a reflection $f'$ of $f$ in $\Lambda$. 
Thus by Theorem~\ref{theorem_lsalg_11} there exists a unique $\Lambda$-compatible monoid operation $\bullet$. 
Moreover, $\Lambda' = (X,f',x_0)$ is also a minimal $\mathcal{L}_S$-algebra and the reflection $\bullet'$ of 
$\bullet$ is the unique $\Lambda'$-compatible monoid operation.

Let $r : (X,f,x_0) \to (X,f',x_0)$ be the unique morphism, hence $r(x_0) = x_0$ and $f'_s \circ r = r \circ f_s$
for all $s \in S$. We next show that $r : (X,f',x_0) \to (X,f,x_0)$ is also a morphism: Let 
$X_0 = \{ x \in X : \mbox{$f_s(r(x)) = r(f'_s(x))$ for all $s \in S$} \}$. Then 
\[f_s(r(x_0)) = f_s(x_0) = f'_s(x_0) = f'_s(r(x_0)) = r(f_s(x_0)) = r(f'_s(x_0))\]
for all $s \in S$, and so $x_0 \in X_0$. Moreover, $X_0$ is $f$-invariant: If $x_0 \in X_0$ and $t \in S$ then 
for all $s \in S$
\begin{eqnarray*}
f_s(r(f_t(x))) &=& f_s(f'_t(r(x))\\ 
&=& f'_t(f_s(r(x))) = f'_t(r(f'_s(x))) = r(f_t(f'_s(x)))  = r(f'_s(f_t(x)))
\end{eqnarray*}
and so $f_t(x) \in X_0$. Thus $X_0 = X$, since $\Lambda$ is minimal, i.e., $f_s \circ r = r \circ f'_s$ for 
all $s \in S$. Since $r(x_0) = x_0$ this means that $r : (X,f',x_0) \to (X,f,x_0)$ is a morphism.

Now $(X,f',x_0)$ is minimal and by Theorem~\ref{theorem_initial_w1} $(X,f,x_0)$ is unambiguous and hence by 
Lemma~\ref{lemma_initial_w6} $(X,f',x_0)$ is unambiguous. Thus by Theorem~\ref{theorem_initial_w1} $(X,f',x_0)$ 
is initial and so by Lemma~\ref{lemma_initial_w3} $r : (X,f,x_0) \to (X,f',x_0)$ is an isomorphism. In 
particular, $r : X \to X$ is a bijection and it then follows from Lemma~\ref{lemma_initial_w2} that 
$r : (X,f',x_0) \to (X,f,x_0)$ is also an isomorphism. Moreover, by Lemma~\ref{lemma_initial_w1}~(2) 
$r \circ r : (X,f,x_0) \to (X,f,x_0)$ is a morphism and hence $r \circ r = \id_X$, since $\id_X$ is the 
unique such morphism. 

Finally, by Theorem~\ref{theorem_initial_w1} $(X,f,x_0)$ and $(X,f',x_0)$ are both unambiguous and in 
particular the mappings $f_s$ and $f'_s$ are injective for each $s \in S$. Proposition~\ref{prop_lsalg_11} 
therefore implies that the monoid $(X,\bullet,x_0)$ obeys the left and right cancellation laws
(since $(X,\bullet,x_0)$ obeying the right cancellation law is the same as
$(X,\bullet',x_0)$ obeying the left cancellation law).
\eop

\begin{lemma}\label{lemma_initial_w11}
Let $(X,f,x_0)$, $(Y,g,y_0)$ be $\mathcal{L}_S$-algebras. If $(X,f,x_0)$ is minimal then there is at most 
one morphism $\pi : (X,f,x_0) \to (Y,g,y_0)$. If $(Y,g,y_0)$ is minimal then any morphism 
$\pi : (X,f,x_0) \to (Y,g,y_0)$ is surjective. Finally, if $(X,f,x_0)$ is minimal then a morphism 
$\pi : (X,f,x_0) \to (Y,g,y_0)$ is surjective if and only if $(Y,g,y_0)$ is minimal.
\end{lemma}

\proof 
As already noted, the first statement holds since if $\pi_1,\,\pi_2$ are morphisms from $(X,f,x_0)$ to 
$(Y,g,y_0)$ then the set $X_0 = \{ x \in X : \pi_1(x) = \pi_2(x) \}$ contains $x_0$ and is $f$-invariant. Thus 
$X_0 = X$, since $(X,f,x_0)$ is minimal, i.e., $\pi_1 = \pi_2$. The second statement follows from the fact that 
$\pi(X)$ contains $y_0$ and is $g$-invariant. (If $y = \pi(x) \in \pi(X)$ and $s \in S$ then 
$g_s(y) = g_s(\pi(x)) = \pi(f_s(x)) \in \pi(X)$.) Hence $\pi(X) = Y$, since $(Y,g,y_0)$ is minimal. It remains 
to show that if $(X,f,x_0)$ is minimal and $\pi : (X,f,x_0) \to (Y,g,y_0)$ is surjective then $(Y,g,y_0)$ is 
minimal. Thus consider a $g$-invariant subset $Y_0$ of $Y$ containing $y_0$. Then $\pi^{-1}(Y_0)$ contains $x_0$
and it is $f$-invariant: If $x \in \pi^{-1}(Y_0)$ (which means $\pi(x) \in Y_0$) and $s \in S$ then
$\pi(f_s(x)) = g_s(\pi(x)) \in Y_0$, since $Y_0$ is $g$-invariant, and so $f_s(x) \in \pi^{-1}(Y_0)$. Therefore 
$\pi^{-1}(Y_0) = X$, since $(X,f,x_0)$ is minimal, which implies $Y_0 = Y$, since $\pi$ is surjective. This shows 
that $(Y,g,y_0)$ is minimal.
\eop

Let us call an $\mathcal{L}_S$-algebra $\Lambda = (X,f,x_0)$ \definition{regular} if it is minimal and there 
exists a (unique) $\Lambda$-compatible monoid operation $\bullet$. In this case $(X,\bullet,x_0)$ will be 
referred to as the associated monoid and $\Lambda' = (X,f',x_0)$ (with $f'$ the reflection of $f$ in $\Lambda$) 
as the reflected  $\mathcal{L}_S$-algebra.

\begin{proposition}\label{prop_initial_w1}
Let $(X,f,x_0)$ and $(Y,g,y_0)$ be regular $\mathcal{L}_S$-algebras and suppose there exists a morphism 
$\pi : (X,f,x_0) \to (Y,g,y_0)$. Then $\pi$ is also a morphism $\pi : (X,f',x_0) \to (Y,g',y_0)$ of the 
reflected $\mathcal{L}_S$-algebras and a homomorphism $\pi : (X,\bullet,x_0) \to (Y,\diamond,y_0)$ of the 
associated monoids.
\end{proposition}

\proof 
We first show $\pi$ is a homomorphism of the associated monoids, and for this consider the set
$X_0 = \{ x \in X : \mbox{$\pi(x \bullet x') = \pi(x) \diamond \pi(x')$ for all $x' \in X$} \}$, which contains 
$x_0$, since $\pi(x_0 \bullet x') = \pi(x') = y_0 \diamond \pi(x') = \pi(x_0) \diamond \pi(x')$ for all 
$x' \in X$. Moreover, $X_0$ is $f$-invariant: If $x \in X_0$ and $s \in S$ then
\begin{eqnarray*}
\pi(f_s(x) \bullet x') &=& \pi(f_s(x \bullet x')) = g_s(\pi(x \bullet x'))\\
&=& g_s(\pi(x) \diamond \pi(x')) = g_s(\pi(x)) \diamond \pi(x') = \pi(f_s(x)) \diamond \pi(x') 
\end{eqnarray*}
for all $x' \in X$, and so $f_s(x) \in X_0$. Thus $X_0$, since $(X,f,x_0)$ is minimal, which shows that $\pi$ 
is a homomorphism. Now if $x \in X$ and $s \in S$ then by Proposition~\ref{prop_lsalg_11}
\[ g'_s(\pi(x)) = \pi(x) \diamond g_s(y_0) 
= \pi(x) \diamond \pi(f_s(x_0)) = \pi(x \bullet f_s(x_0)) = \pi(f'_s(x))
\] 
and hence $g'_s \circ \pi = \pi \circ f'_s$ for all $s \in S$. This shows that $\pi$ is a morphism of the 
reflected $\mathcal{L}_S$-algebras.
\eop

\begin{lemma}\label{lemma_initial_w12}
Let $(X,f,x_0)$ be a regular and $(Y,g,y_0)$ a minimal $\mathcal{L}_S$-algebra and suppose there exists a 
morphism $\pi : (X,f,x_0) \to (Y,g,y_0)$. Then:
\begin{evlist}{8pt}{6pt}
\item[(1)]
$\pi(f_s(x_1)) = \pi(f_s(x_2))$ holds for all  $s \in S$ whenever $\pi(x_1) = \pi(x_2)$. 

\item[(2)]
$\pi(x \bullet x_1) = \pi(x \bullet x_2)$ holds for all $x \in X$ whenever $\pi(x_1) = \pi(x_2)$. 
\end{evlist}
\end{lemma}

\proof 
(1)\enskip
If $\pi(x_1) = \pi(x_2)$ then $\pi(f_s(x_1)) = g_s(\pi(x_1)) = g_s(\pi(x_2)) = \pi(f_s(x_2))$.

(2)\enskip
Let $x_1,\,x_2 \in X$ with $\pi(x_1) = \pi(x_2)$ and consider the set $X_0$ of those $x \in X$ for which 
$\pi(x \bullet x_1) = \pi(x \bullet x_2)$. Then 
$\pi(x_0 \bullet x_1) = \pi(x_1) = \pi(x_2) = \pi(x_0 \bullet x_2)$ and so $x_0 \in X$. Moreover, $X_0$ is 
$f$-invariant: If $x \in X_0$ and $s \in S$ then
\begin{eqnarray*}
\pi(f_s(x) \bullet x_1) &=& \pi(f_s(x \bullet x_1))\\ 
&=& g_s(\pi(x \bullet x_1)) = g_s(\pi(x \bullet x_2)) = \pi(f_s(x \bullet x_2)) = \pi(f_s(x) \bullet x_2)
\end{eqnarray*}
and so $f_s(x) \in X_0$. Thus $X_0 = X$, since $(X,f,x_0)$ is minimal.
\eop

\begin{proposition}\label{prop_initial_w2}
Let $(X,f,x_0)$ be a regular and $(Y,g,y_0)$ a minimal $\mathcal{L}_S$-algebra and suppose that there
exists a morphism $\pi : (X,f,x_0) \to (Y,g,y_0)$. Then the following are equivalent:
\begin{evlist}{8pt}{6pt}
\item[(1)]
$\pi(f'_s(x_1)) = \pi(f'_s(x_2))$ holds for all  $s \in S$ whenever $\pi(x_1) = \pi(x_2)$. 

\item[(2)]
$\pi(x_1 \bullet x) = \pi(x_2 \bullet x)$ holds for all $x \in X$ whenever $\pi(x_1) = \pi(x_2)$. 

\item[(3)]
$(Y,g,y_0)$ is regular.
\end{evlist}
\end{proposition}

\proof 
(2) $\Rightarrow$ (3):\enskip
Let $x_1,\,x'_1,\,x_2,\,x'_2 \in X$ with $\pi(x_1) = \pi(x'_1)$ and $\pi(x_2) = \pi(x'_2)$. Then 
$\pi(x_1\bullet x_2) = \pi(x'_1\bullet x_2)$ and by Lemma~\ref{lemma_initial_w12}~(2)
$\pi(x'_1\bullet x_2) = \pi(x'_1\bullet x'_2)$ and so $\pi(x_1\bullet x_2) = \pi(x'_1\bullet x'_2)$. Thus, since 
by Lemma~\ref{lemma_initial_w11} $\pi$ is surjective, there exists a unique binary operation $\diamond$ on 
$Y$ such that $\pi(x_1 \bullet x_2) = \pi(x_1) \diamond \pi(x_2)$ for all $x_1,\,x_2 \in X$. Then $\diamond$ 
is a monoid operation on $(Y,y_0)$: Let $y_1,\,y_2,\,y_3 \in Y$ and choose $x_1,\,x_2,\,x_3 \in X$ with
$\pi(x_j) = y_j$ for $j = 1,\,2,\,3$. It follows that
\begin{eqnarray*}
(y_1 \diamond y_2) \diamond y_3 &=& (\pi(x_1) \diamond \pi(x_2)) \diamond \pi(x_3)\\
 &=&  \pi(x_1 \bullet x_2) \diamond \pi(x_3) = \pi((x_1 \bullet x_2) \bullet x_3)
= \pi(x_1 \bullet (x_2 \bullet x_3))\\  
&=&  \pi(x_1) \diamond \pi(x_2 \bullet x_3) =  \pi(x_1) \diamond (\pi(x_2) \diamond \pi(x_3)) 
= y_1 \diamond (y_2 \diamond y_3) 
\end{eqnarray*}
and so $\diamond$ is associative. Also, if $y = \pi(x)$ then
\[y \diamond y_0 = \pi(x) \diamond \pi(x_0) = \pi(x \bullet x_0) = \pi(x) = y\]
and in the same way $y_0 \diamond y = y$. Finally, if $y = \pi(x)$ and $s \in S$ then
\begin{eqnarray*}
g_s(y) &=& g_s(\pi(x)) = \pi(f_s(x))\\ 
&=& \pi(f_s(x_0) \bullet x) = \pi(f_s(x_0)) \diamond \pi(x) = g_s(\pi(x_0)) \diamond y = g_s(y_0) \diamond y \;,
\end{eqnarray*}
which shows that $g_s$ is a translation in $(Y,\diamond,y_0)$ for each $s \in S$. Hence $(Y,g,y_0)$ is regular.

(3) $\Rightarrow$ (1):\enskip
By Proposition~\ref{prop_initial_w1} $\pi$ is also a morphism of the reflected $\mathcal{L}_S$-algebras 
$\pi : (X,f',x_0) \to (Y,g',y_0)$. Thus if $x_1,\,x_2 \in X$ with $\pi(x_1) = \pi(x_2)$ then 
\[ \pi(f'_s(x_1)) =  g'_s(\pi(x_1)) = g'_s(\pi(x_2)) = \pi(f'_s(x_2))\;.\]

(1) $\Rightarrow$ (2):\enskip
Consider the set $X_0$ of those $x \in X$ for which $\pi(x_1 \bullet x) = \pi(x_2 \bullet x)$ holds whenever 
$x_1,\,x_2 \in X$ with $\pi(x_1) = \pi(x_2)$, and so in particular $x_0 \in X_0$. Let $x \in X_0$  and 
$x_1,\,x_2 \in X$ with $\pi(x_1) = \pi(x_2)$. Then $\pi(x_1 \bullet x) = \pi(x_2 \bullet x)$, since 
$x \in X_0$, and so by (1) $\pi(f'_s(x_1 \bullet x)) = \pi(f'_s(x_2 \bullet x))$. But 
$f'_s(x' \bullet x) = x' \bullet f'_x(x)$ for all $x' \in X$ and hence
$\pi(x_1 \bullet f'_s(x)) = \pi(x_2 \bullet f'_s(x))$, which shows that $f'_s(x) \in X_0$, i.e., $X_0$ is 
$f'$-invariant. Thus $X_0 = X$, since $(X,f',x_0)$ is minimal, which implies that (2) holds.
\eop

Let $(X,f,x_0)$ be an $\mathcal{L}_S$-algebra, let $(Y,y_0)$ be a pointed set and $p : X \to Y$ be a surjective 
mapping with $p(x_0) = y_0$. Then there is an equivalence relation $\approx$ on $X$ with $x_1 \approx x_2$ if 
and only if $p(x_1) = p(x_2)$. Conversely, if we start with an equivalence relation $\approx$ on $X$, let $Y$ 
be the set of equivalence classes and $p$ be the mapping which assigns to each element $x$ the equivalence class 
$[x]$ to which it belongs then $p : X \to Y$ is surjective with $p(x_0) = y_0$, where $y_0 = [x_0]$. Now there 
exists a mapping $g : S \times Y \to Y$ so that $p : (X,f,x_0) \to (Y,g,y_0)$ is a morphism if and only if 
$p(f_s(x_1)) = p(f_s(x_2))$ for all $s \in S$ whenever $p(x_1) = p(x_2)$, or, what is equivalent, if and only if
$f_s(x_1) \approx f_s(x_2)$  for all $s \in S$ whenever $x_1 \approx x_2$. 

Suppose this requirement is met. If $(X,f,x_0)$ is minimal then by Lemma~\ref{lemma_initial_w11} $(Y,g,y_0)$ 
is also minimal. Moreover, Proposition~\ref{prop_initial_w2} implies that if $(X,f,x_0)$ is regular then 
$(Y,g,y_0)$ is regular if and only if $p(f'_s(x_1)) = p(f'_s(x_2))$ for all $s \in S$ whenever 
$p(x_1) = p(x_2)$ (or, what is the same, if and only if  $f'_s(x_1) \approx f'_s(x_2)$ for all $s \in S$ 
whenever $x_1 \approx x_2$). In this case it follows from Proposition~\ref{prop_initial_w1} that
$\pi : (X,\bullet,x_0) \to (Y,\diamond,y_0)$ is a homomorphism of the associated monoids. 

We end the section by considering a class of $\mathcal{L}_S$-algebras defined in terms of suitable subsets
of a given initial $\mathcal{L}_S$-algebra.
In what follows let $(X,f,x_0)$ be an initial $\mathcal{L}_S$-algebra (and so by 
Theorem~\ref{theorem_initial_w1} $(X,f,x_0)$ is minimal and unambiguous) and let $A$ be a subset of $X$ 
containing $x_0$ and such that $X \setminus A$ is $f$-invariant. Put
\[ \partial A = \{ x \in X \setminus A\,:\,\mbox{$x = f_s(x')$ for some $x' \in A$ and some $s \in S$}\} \]
and let $\bar{A} = A \cup \partial A$. 

\begin{lemma}\label{lemma_initial_w13}
$f_s(A) \subset \bar{A}$ and $f_s(X \setminus A) \subset X \setminus \bar{A}$ for each $s \in S$.
In particular, the set $X \setminus \bar{A}$ is also $f$-invariant.
\end{lemma}

\proof
The first statement (that $f_s(A) \subset \bar{A}$) follows from the definition of $\partial A$. Now if 
$f_s(x) \in \partial A$ for some $s \in S$ then $x \in A$ (since $f_s(x) = f_t(x')$ for some $x' \in A$, 
$t \in S$ and then $s = t$ and $x = x'$ because $(X,f,x_0)$ is unambiguous). Hence if $x \in X \setminus A$ 
and $s \in S$ then $f_s(x) \notin \partial A$, and also $f_s(x) \in X \setminus A$, since $X \setminus A$ is 
$f$-invariant. Thus $f_s(x) \in X \setminus \bar{A}$ for all $x \in X \setminus A$. 
\eop

Since $f_s(A) \subset \bar{A}$ we can define for each $s \in S$ a mapping $f^A_s : \bar{A} \to \bar{A}$ by
\[   f^A_s(x) = \left\{ \begin{array}{cl}
                        f_s(x)  &\ \mbox{if $x \in A$}\;,\\
                          x     &\ \mbox{if $x \in \partial A$}\;.
                 \end{array} \right. \]
This gives us an $\mathcal{L}_S$-algebra $(\bar{A},f^A,x_0)$.

\begin{lemma}\label{lemma_initial_w14}
The $\mathcal{L}_S$-algebra $(\bar{A},f^A,x_0)$ is minimal.
\end{lemma}

\proof 
Let $A_0$ be an $f^A$-invariant subset of $\bar{A}$ containing $x_0$, and consider the subset
$X_0 = A_0 \cup (X \setminus \bar{A})$ of $X$. Then $x_0 \in X_0$ and $X_0$ is $f$-invariant: Let $x \in X_0$ 
and $s \in S$. If $x \in X \setminus A$ then by Lemma~\ref{lemma_initial_w13} 
$f_s(x) \in X \setminus \bar{A} \subset X_0$; on the other hand, if $x \in A$ then $x \in A_0$ and so
$f_s(x) = f^A_s(x) \in X_0$. Thus in both cases $f_s(x) \in X_0$. Hence $X_0 = X$, since $(X,f,x_0)$ is minimal, 
which implies that $A_0 = \bar{A}$. This shows that $(\bar{A},f^A,x_0)$ is minimal. \eop

\begin{theorem}\label{theorem_initial_w3}
The $\mathcal{L}_S$-algebra $(\bar{A},f^A,x_0)$ regular if and only if the set $X \setminus A$ is also 
$f'$-invariant (with $f'$ the reflection of $f$ in $(X,f,x_0)$).
\end{theorem}

The proof of Theorem~\ref{theorem_initial_w3} requires some preparation. Since $(X,f,x_0)$ is initial there 
exists a unique morphism $p : (X,f,x_0) \to (\bar{A},f^A,x_0)$. Hence $p(x_0) = x_0$ and 
$p \circ f_s = f^A_s \circ p$, i.e., for all $s \in S$
\[   p(f_s(x)) = \left\{ \begin{array}{cl}
                        f_s(p(x))  &\ \mbox{if $p(x) \in A$}\;,\\
                          p(x)     &\ \mbox{if $p(x) \in \partial A$\;.}
                 \end{array} \right. \]

\begin{lemma}\label{lemma_initial_w15}
(1)\enskip 
$p(x) = x$ for all $x \in \bar{A}$.

(2)\enskip
$p(x) \in \partial A$ for all $x \in X \setminus A$.

(3)\enskip
$p(f_s(x)) = p(x)$ for all $x \in X \setminus A$, $s \in S$.

\end{lemma}

\proof 
(1)\enskip
Let $X_0 = \{ x \in X : p(x) = x \} \cup (X \setminus A)$. Then $x_0 \in X_0$, since $p(x_0) = x_0$,
and $X_0$ is $f$-invariant: Let $x \in X_0$ and $s \in S$; if $x \in X \setminus A$ then 
$f_s(x) \in X \setminus A \subset X_0$, since $X \setminus A$ is $f$-invariant. On the other hand, if 
$x \in A$ then $p(x) = x$ (and so in particular $p(x) \in A$) and hence $p(f_s(x)) = f_s(p(x)) = f_s(x)$, 
which again means $f_s(x) \in X_0$. Thus $X_0 = X$, since $(X,f,x_0)$ is minimal, and this shows that $p(x) = x$ 
for all $x \in A$. Moreover, if $x \in \partial A$ then $x = f_s(x')$ for some $x' \in A$ and some $s \in S$ and 
then, since $p(x') = x'$, it follows that $p(x) = p(f_s(x')) = f_s(p(x')) = f_s(x') = X$, i.e., $p(x) = x$ holds 
for all $x \in \partial A$.
 
(2)\enskip
Let $X_0 = \{ x \in X : p(x) \in \partial A \} \cup A$. Then $x_0 \in X_0$, since $x_0 \in A$, and $X_0$ is 
$f$-invariant: Let $x \in X_0$ and $s \in S$; if $x \in X \setminus A$ then $p(x) \in \partial A$ and so 
$p(f_s(x)) = p(x) \in \partial A$, i.e., $f_s(x) \in X_0$. On the other hand, if $x \in A$ then either 
$f_s(x) \in A$, in which case $f_s(x) \in X_0$, or $f_s(x) \in \partial A$, in which case by (1) 
$p(f_s(x)) = f_s(x) \in \partial A$, and again $f_s(x) \in X_0$. Thus $X_0 = X$, since $(X,f,x_0)$ is minimal, 
and this shows that $p(x) \in \partial A$ for all $x \in X \setminus A$. 

(3)\enskip
This follows from (2), since $p(f_s(x)) = p(x)$ whenever $p(x) \in \partial A$.
\eop

By Lemma~\ref{lemma_initial_w15}~(1) the mapping $p : X \to \bar{A}$ is surjective and so
Lemma~\ref{lemma_initial_w11} confirms that $(\bar{A},f^A,x_0)$ is a minimal $\mathcal{L}_S$-algebra.

\begin{lemma}\label{lemma_initial_w16}
(1)\enskip
For each $x \in X \setminus A$ there exists $x' \in X$ with $x = x' \bullet p(x)$.

(2)\enskip 
If $x \in X \setminus A$ then $x' \bullet x \in X \setminus A$ and $p(x' \bullet x) = p(x)$ for all $x' \in X$.
\end{lemma}

\proof
(1)\enskip 
Let $X_0 = A \cup \{ x \in X \setminus A : \mbox{there exists $x' \in X$ with $x = x' \bullet p(x)$} \}$. Then 
$x_0 \in A \subset X_0$, and $X_0$ is $f$-invariant: Consider $x \in X_0$ and $s \in S$ and we can assume that 
$f_s(x) \in X \setminus A$. If $x \in A$ then $f_s(x) \in \partial A$ and so by Lemma~\ref{lemma_initial_w15}~(1)
$f_s(x) = p(f_s(x)) = x_0 \bullet p(f_s(x))$, i.e., $f_s(x) \in X_0$. On the other hand, if $x \in X \setminus A$ 
then $x = x' \bullet p(x)$ for some $x' \in X$ and by Lemma~\ref{lemma_initial_w15}~(3) $p(x) = p(f_s(x))$ and 
hence $f_s(x) = f_s(x' \bullet p(x)) = f_s(x') \bullet p(x) = f_s(x') \bullet p(f_s(x))$, and again
$f_s(x) \in X_0$. Thus $X_0 = X$, since $(X,f,x_0)$ is minimal, and this shows that 
for each $x \in X \setminus A$ there exists $x' \in X$ with $x = x' \bullet p(x)$.

(2)\enskip
Fix $x \in X \setminus A$ and let 
$X_0 = \{ x' \in X : \mbox{$x' \bullet x \in X \setminus A$ and $p(x' \bullet x) = p(x)$} \}$. Then
$X_0$ contains $x_0$, since $x_0 \bullet x = x$, and it is $f$-invariant: If $x' \in X_0$ and $s \in S$ then
$f_s(x') \bullet x = f_s(x' \bullet x) \in X \setminus A$, since $X \setminus A$ is $f$-invariant, and by 
Lemma~\ref{lemma_initial_w15}~(3) $p(f_s(x') \bullet x) = p(f_s(x' \bullet x)) = p(x' \bullet x) = p(x)$, and so 
$f_s(x') \in X_0$. Thus $X_0 = X$, since $(X,f,x_0)$ is minimal, and this shows that 
$x' \bullet x \in X \setminus A$ and $p(x' \bullet x) = p(x)$ for all $x' \in X$.
\eop

\begin{lemma}\label{lemma_initial_w17}
$(\bar{A},f^A,x_0)$ is regular if and only if $p(f'_s(x)) = p(f'_s(p(x))$ for all $x \in X \setminus A$, 
$s \in S$.
\end{lemma}

\proof 
By Proposition~\ref{prop_initial_w2} the $\mathcal{L}_S$-algebra $(\bar{A},f^A,x_0)$ is regular if and only if
$p(f'_s(x_1)) = p(f'_s(x_2))$ for all $s \in S$ whenever $x_1,\,x_2 \in X$ with $p(x_1) = p(x_2)$. Thus if 
$(\bar{A},f^A,x_0)$ is regular then $p(f'_s(x)) = p(f'_s(p(x))$ for all $x \in X \setminus A$, $s \in S$, since 
by Lemma~\ref{lemma_initial_w15} (1) and (2) $p(x) = p(p(x))$ for all $x \in X \setminus A$. Suppose conversely 
$p(f'_s(x)) = p(f'_s(p(x))$ for all $x \in X \setminus A$, $s \in S$, and consider $x_1,\,x_2 \in X$ with 
$p(x_1) = p(x_2)$. We will show that $p(f'_s(x_1)) = p(f'_s(x_2))$. Now if at least one of $x_1$ and $x_2$ is in 
$A$ then by Lemma~\ref{lemma_initial_w15} $p(x_1) = p(x_2)$ holds only if $x_1 = x_2$ and in this case 
$p(f'_s(x_1)) = p(f'_s(x_2))$ holds trivially. We can therefore assume that $x_1,\,x_2 \in X \setminus A$ and 
hence $p(f'_s(x_1)) = p(f'_s(p(x_1)) = p(f'_s(p(x_2)) = p(f'_s(x_2))$ for all $s \in S$. \eop

\textit{Proof of Theorem~\ref{theorem_initial_w3}:\ }
Suppose first that $X \setminus A$ is  $f'$-invariant. Let $x \in X \setminus A$ and put $\hat{x} = p(x)$.
By Lemma~\ref{lemma_initial_w15}~(2) $\hat{x} \in \partial A \subset X \setminus A$ and by 
Lemma~\ref{lemma_initial_w16}~(1) there exists $x' \in X$ with $x = x' \bullet \hat{x}$. It follows that
$f'_s(x) = f'_s(x' \bullet \hat{x}) = x' \bullet f'_s(\hat{x})$ and $f'_s(\hat{x}) \in X \setminus A$, since 
$X \setminus A$ is $f'$-invariant; hence by Lemma~\ref{lemma_initial_w16}~(2) 
\[ p(f'_s(x)) = p(x' \bullet f'_s(\hat{x})) = p(f'_s(\hat{x})) = p(f'_s(p(x)))\] 
for all $s \in S$. Therefore by Lemma~\ref{lemma_initial_w17} $(\bar{A},f^A,x_0)$ is regular.

Suppose conversely that $(\bar{A},f^A,x_0)$ is regular and let $x \in X \setminus A$ and $s \in S$. By 
Lemma~\ref{lemma_initial_w17} $p(f'_s(x)) = p(f'_s(p(x))$ and thus if $f'_s(x) \in A$ then by
Lemma~\ref{lemma_initial_w15} $f'_s(x) = f'_s(p(x))$, which implies that $x = p(x)$, since by
Theorems \ref{theorem_initial_w1} and \ref{theorem_initial_w2} $f'$ is unambiguous, which in turn implies that 
$x \in \partial A$. This shows $f'_s(x) \in X \setminus A$ for all $x \in X \setminus \bar{A}$ and so it remains 
to show that $f'_s(x) \in X  \setminus A$ whenever $x \in \partial A$. Let $x' \in X \setminus \{x_0\}$; 
by Lemma~\ref{lemma_initial_w16} $x' \bullet x \in X \setminus A$ and $p(x' \bullet x) = p(x) = x$ and hence 
$p(f'_s(x)) = p(f'_s(p(x' \bullet x))) = p(f'_s(x' \bullet x))$. If $f'_s(x) \in A$ then, as above, it would 
follow that $f'_s(x) = f'_s(x' \bullet x)$ and therefore $x_0 \bullet x = x = x' \bullet x$, again since $f'$ is 
unambiguous. But this is not possible, since by Theorem~\ref{theorem_initial_w2} the right cancellation law holds 
in $(X,\bullet,x_0)$. Thus $f'_s(x) \in X \setminus A$.
\eop

\bigskip
\bigskip

{\sc Fakult\"at f\"ur Mathematik, Universit\"at Bielefeld}\\
{\sc Postfach 100131, 33501 Bielefeld, Germany}\\
\textit{E-mail address:} \texttt{preston@math.uni-bielefeld.de}\\
